\documentclass[11pt]{article}
\input xy
\xyoption{all}
\usepackage{amsmath, amsthm, amssymb}

\def\ddefi#1{

{\noindent{\bf Definition #1.}

}}
\def\nota#1#2{

{\noindent{{\bf Notation #1}#2.}

}}
\def\lem#1{

{\noindent{\bf Lemma #1.}

}}
\def\prop#1{

{\noindent{\bf Proposition #1.}

} }
\def\teo#1{

{\noindent{\bf Theorem #1.}

}}
\def\cor#1{

{\noindent{\bf Corollary #1.}

}}
\def\dim{

{\noindent{\bf Proof:} }}

\def\sqr#1#2{{\vcenter{\vbox{\hrule height .#2pt
                             \hbox{\vrule width .#2pt height#1pt\kern#1pt
                                   \vrule width .#2pt}
                             \hrule height .#2pt}}}}

\def\rem#1#2{

{\noindent{\bf Remark #1}#2.

}}

\def\u{{\cal U}}
\def\uq{{\cal U}_q}
\def\udr{{\cal U}_q^{Dr}}
\def\udj{{\cal U}_q^{DJ}}
\def\a{{\cal A}}
\def\d{\delta}
\def\ua{{\cal U}_{\a}^{DJ}}
\def\uad{{\cal U}_{\a}^{Dr}}
\def\uado{{\cal U}_{\a}^{Dr,0}}
\def\udjp{{\cal U}_{\a}^{DJ,+}}
\def\f{{\cal F}}

\def\iz{I_{\Z}}

\def\i{{\cal I}}
\def\v{{\cal{V}}}
\def\sy{{\cal S}}
\def\gothg{\frak g}
\def\gothh{\frak h}
\def\gothl{\frak l}

\def\n{\noindent}

\def\compo{\,{\scriptstyle\circ}\,}

\def\N{{\Bbb{N}}
}
\def\Z{{\Bbb{Z}}
}

\def\C{{\Bbb{C}}
}

\def\intr{\S 0}

\def\specgen{0.1}

\def\gnrntt{\S 1}

\def\gennu{1.1}
\def\rsy{1.2}
\def\gent{{\bf 1.3}}
\def\genno{1.4}
\def\nbr{1.5}

\def\preldj{\S 2}

\def\QNTALG{2.1}
\def\qremu{2.2}
\def\grgr{2.3}
\def\qremd{2.4}

\def\preldr{\S 3}

\def\drdef{3.1}
\def\drnot{3.2}
\def\drrem{3.3}
\def\drom{3.4}
\def\drromrem{3.5}
\def\uqz{3.6}
\def\dremrelu{3.7}
\def\stagg{3.8}
\def\drnotu{3.9}
\def\dremd{3.10}

\def\prelpsi{\S 4}

\def\defpsi{4.1}
\def\repsi{4.2}
\def\copsi{4.3}

\def\trngfindim{\S 5}

\def\prpsid{5.1}
\def\prpsit{5.2}
\def\repsit{5.3}
\def\epp{5.4}
\def\pip{5.5}

\def\ntgf{\S 6}

\def\spnot{6.1}
\def\firem{6.2}
\def\insint{6.3}
\def\pfp{6.4}
\def\imredu{6.5}
\def\uffnot{6.6}
\def\pshom{6.7}

\def\spquno{\S 7}

\def\imredd{7.1}
\def\imredt{7.2}
\def\imredq{7.3}
\def\aggde{7.4}
\def\aggre{7.5}
\def\imredc{7.6}
\def\imreds{7.7}
\def\epip{7.8}
\def\imredo{7.9}
\def\imleu{7.10}
\def\imcou{7.11}
\def\imled{7.12}
\def\imredn{7.13}
\def\emo{7.14}
\def\emu{7.15}
\def\emd{7.16}
\def\nqud{7.17}
\def\ntnz{7.18}
\def\relsem{7.19}
\def\nrku{7.20}
\def\nlmu{7.21}
\def\ntzu{7.22}
\def\nrkt{7.23}
\def\nrkq{7.24}
\def\ntzd{7.25}
\def\nrkc{7.26}
\def\imlec{7.27}
\def\nrks{7.28}
\def\lpgen{7.29}

\def\affkm{\S 8}

\def\ddd{8.1}
\def\loab{8.2}
\def\lopm{8.3}
\def\lgemb{8.4}
\def\lach{8.5}
\def\pmstr{8.6}
\def\lpacg{8.7}
\def\efhom{8.8}
\def\dvfi{8.9}
\def\lacf{8.10}
\def\lpacef{8.11}
\def\gtg{8.12}
\def\eifi{8.13}
\def\impser{8.14}
\def\lrgm{8.15}
\def\genw{8.16}
\def\gnaut{8.17}
\def\efexp{8.18}
\def\porb{8.19}
\def\ult{8.20}
\def\rdinsq{8.21}

\def\cncl{\S 9}

\def\tttu{9.1}
\def\tttd{9.2}
\def\tttq{9.3}
\def\tttc{9.4}
\def\ttto{9.5}
\def\tttt{9.6}

\def\matsu{\cite{13}}

\def\beck{\cite{2}}
\def\bbk{\cite{3}}
\def\realdrel{\cite{4}}
\def\damcina{\cite{5}}
\def\drld{\cite{6}}
\def\drr{\cite{7}}
\def\hern{\cite{8}}
\def\jm{\cite{9}}
\def\kac{\cite{10}}
\def\levsb{\cite{11}}
\def\lusz{\cite{12}}

\begin{document} 
\title{FROM THE DRINFELD REALIZATION TO THE DRINFELD-JIMBO PRESENTATION OF 
AFFINE QUANTUM
ALGEBRAS: THE INJECTIVITY.}
\date{}
\author{Ilaria Damiani}
\maketitle

\begin{abstract} \noindent In this paper the surjective homomorphism $\psi$ (see \realdrel) from the Drinfeld realization $\udr$ to the Drinfeld and Jimbo presentation $\udj$ of affine quantum algebras is proved to be injective.

\n A consequence of the arguments used in the paper is the triangular decomposition of the Drinfeld realization of affine quantum algebras also in the twisted case.

\n A presentation of the affine Kac-Moody algebras in terms of the ``Drinfeld generators'' is also provided.
\end{abstract}

\vskip .5 truecm
\noindent {\bf{\intr.\ {\bf INTRODUCTION.}}}
\vskip .5truecm

\noindent Let $X_{\tilde n}^{(k)}$ be a Dynkin diagram of affine type, $\udj=\udj(X_{\tilde n}^{(k)})$ the quantum algebra introduced by Drinfeld and Jimbo (see \drr\ and \jm), $\udr=\udr(X_{\tilde n}^{(k)})$ its Drinfeld realization (see \drld).

\noindent This paper concludes the proof that $\udj$ and $\udr$ are isomorphic. More precisely, in \realdrel\ a homomorphism $\psi:\udr\to\udj$ was defined (following \beck\ for the untwisted case), and proved to be surjective; previous attempts to give a complete proof that these two algebras are isomorphic are also discussed in  \realdrel. Here is a proof of the injectivity of $\psi$.

\noindent As in \beck, the idea of the proof is recovering the injectivity of $\psi$ from that of its specialization at 1, based on the following:
\vskip .3 truecm
\prop{\specgen}
\noindent Let $\a=\C[q]_{(q-1)}$ be the localization of $\C[q]$ at $(q-1)$, $M$ a finitely generated $\a$-module, $N$ a free $\a$-module, $f:M\to N$ a homomorphism of $\a$-modules, $f_1:M/(q-1)M\to N/(q-1)N$ the $\a/(q-1)=\C$-linear homomorphism induced by $f$.

If $f_1$ is injective, so is $f$.
\dim
\noindent $\a$ is a local principal ideal domain; $f(M)$ is a finitely generated $\a$-submodule of $N$, hence a free $\a$-module, so that 
there exists $g:f(M)\to M$ such that $f\compo g=id_{f(M)}$.

\noindent Of course $ker(f)$ is a finitely generated $\a$-module, $M=ker(f)\oplus Im(g)
$, $ker(f)/(q-1)ker(f)\hookrightarrow M/(q-1)M$ and $ker(f)/(q-1)ker(f)\subseteq ker(f_1)=\{0\}$.

\noindent Then $(q-1)ker(f)=ker(f)$, so that $ker(f)=\{0\}$ (Nakayama lemma).

Remark that the hypothesis that $M$ is finitely generated over $\a$ is necessary, as it can be seen from the simple counterexample $f:\C(q)\to\{0\}$.
\vskip .3 truecm
\noindent  The problem faced in the present paper is reducing to a situation where this argument works. 

\n Consider the (well defined) commutative diagram 
$$\xymatrix{
&\f_+/\i_+\ar[r]^{f}\ar[d]_{\tilde\psi}&\udr\ar[d]^{\psi}\\&\udjp\ar@{^{(}->}[r]&\udj&{}}$$
where $\udj$ and $\udr$ are respectively the Drinfeld-Jimbo presentation and the Drinfeld realization of a quantum affine algebra (see sections \preldj\ and \preldr), $\udjp$ is the integer form of the positive part of $\udj$ (remark \qremd), $\f_+$ is the free $\a$-algebra generated by $\{X_{i,r}^+|i\in I_0, r\geq 0\}$ and $\i_+$ is the ideal of $\f_+$ generated by the relations $(ZX_+^+,DR_+, S_+,U3_+)$ (see notations \drnotu\ and \spnot).

\n 
The plan of the proof is showing that the injectivity of $\tilde\psi$ implies the injectivity of $\psi$ (see proposition \prpsit\ and corollary \pfp, ii)) and at the same time that the conditions of proposition \specgen\ hold for the homogeneous components of  $\tilde\psi:\f_+/\i_+\to\udjp$ (see remark \pshom), so  that $
\psi$ is injective if $\tilde\psi_1$ - the specialization at 1  of $\tilde\psi$ - is injective. This turns our problem into the study of $\tilde\psi_1$, which is found out to be injective through a careful analysis of the classical (non quantum) affine Kac-Moody case (see remark \imredc\ and corollary \rdinsq).

A) It is well known that $\udjp$ is a free $\a$-module (see remark \qremd,ii)); it is straightforward to see that $\f_+/\i_+=\oplus_{\alpha\in Q_+}(\f_+/\i_+)_{\alpha}$ where each $(\f_+/\i_+)_{\alpha}$ is a finitely generated $\a$-module (see remark \firem, ii)); finally $\tilde\psi$ is trivially $Q$-homogeneous: then proposition \specgen\ applies and $\tilde\psi$ {\it is injective if $\tilde\psi_1$ is injective.}

B) Of course $\f_+/\i_+$, but also $\udjp$, can be easily described through a {\it presentation by generators and relations} (it is well known that $\udjp$ is generated by $\{E_i|i\in I\}$ with relations $(SE)$, see remark \qremd, iv)). Then their specializations at 1 are also immediate to describe by generators and relations (see remarks \imredd\ and \aggre), and $\tilde\psi_1$ is explicitly known on the generators.
\n Section \affkm\ is devoted to prove that $\tilde\psi_1$ {\it is injective}. Since the specialization at 1 of $\udjp$ is well know (it is the enveloping algebra of the positive part of the Kac-Moody algebra), the proof consists in the study of the classical (non quantum) situation, through a careful analysis of the specialization at 1 of $\f_+/\i_+$ (see corollary \lpgen\ and section \affkm). In particular this analysis leads also to a 
``Drinfeld realization'' of the affine Kac-Moody algebras (see theorem \tttt). 

C) On the other hand $f(\f_+/\i_+)$ generates over $\C(q)$ a subalgebra $\uq^{Dr,+,+}$ of $\uq^{Dr,+}\subseteq\udr$; since $f(\f_+/\i_+)$ is direct sum of finitely generated $\a$-modules, it is an integer form of $\uq^{Dr,+,+}$ (see remark \insint). So the injectivity of $\tilde\psi$ implies that  $\psi\big|_{\uq^{Dr,+,+}}$ {\it is injective} (see corollary \pfp, ii)).

\n But the injectivity of $\tilde\psi$ (then of $f$) implies also that $\f_+/\i_+\cong f(\f_+/\i_+)$, that is it provides a presentation by generators and relations of the integer form of $\uq^{Dr,+,+}$ (see corollary \pfp, i)).

D) Why does the injectivity of $\psi\big|_{\uq^{Dr,+,+}}$ imply the injectivity of $\psi$?

\n To answer this question we study the connection between the PBW basis of $\udj$ and the tensor product $\uq^{Dr,-,-}\otimes\uq^{Dr,0}\otimes\uq^{Dr,+,+}$ (see proposition \prpsid), recalling that  
$\uq^{Dr,+}$ can be recovered from $\uq^{Dr,+,+}$ by ``translations'' (see remark \drromrem, vii)).

\n With these tools it is easy to conclude finally that the injectivity of $\psi\big|_{\uq^{Dr,+,+}}$ implies the injectivity of  $\psi$. At the same time it implies also the triangular decomposition of $\udr$ (see proposition \prpsit).

\vskip .3 truecm
\n I want to thank Velleda Baldoni for her support and advice and Hiraku Nakajima for his interest and attention.

\vskip .5 truecm
\noindent {\bf{\gnrntt.\ {\bf GENERAL NOTATIONS.}}}
\vskip .5truecm

\n We fix here the general notations that will be used in the paper (for a deeper and more detailed understanding of this setting see \bbk, \kac, \matsu). Further notations are spread out in the next sections, following the exposition.

\vskip .3 truecm

\nota{\gennu}{}
\vskip .3 truecm

\n Following the literature we denote by 

$\hat\gothg=(\gothg\otimes_{\C}\C[t^{\pm 1}])^{\chi}\oplus\C c$ an affine Kac-Moody algebra 

\n with 

Dynkin diagram $\Gamma$ and set of vertices $I=\{0,1,...,n\}$, 

Cartan matrix $A=(a_{ij})_{i,j\in I}$, 


root lattice $Q=\oplus_{i\in I}\Z\alpha_i$ and positive root lattice $Q_+=\oplus_{i\in I}\N\alpha_i$, 

root system (with real and imaginary roots) $\Phi=\Phi^{{\rm{re}}}\cup\Phi^{{\rm{im}}}$, 

root system with multiplicities $\hat\Phi
$,

symmetric bilinear form $(\cdot|\cdot)$ on Q induced by $DA$,    ($D=$diag$(d_i|i\in I)$) with kernel $\Z\delta$ ($\delta\in Q_+$),

Weyl group $W=<s_i:\alpha_j\mapsto\alpha_j-a_{ij}\alpha_i|i\in I>$, 

extended Weyl group $\hat W=W\rtimes {\cal T}$ (${\cal T}\leq Aut(\Gamma)$) with length $l:\hat W\to\N$, 

extended braid group 
with lifting $\hat W\ni w\mapsto
T_w
$, 

\n where 

$\gothg$ is a simple Lie algebra over $\C$ of rank $\tilde n$;

$\chi$ is an automorphism of the Dynkin diagram of $\gothg$ of order $k$;

$A_0=(a_{ij})_{i,j\in I_0}$ ($I=I_0\cup\{0\}\neq I_0$), 

$Q_0=\oplus_{i\in I_0}\Z\alpha_i\subseteq Q$ ($Q_{0,+}=Q_0\cap Q_+$),

$\Phi_0\supseteq\Phi_{0,+} $ and 

$W_0=<s_i|i\in I_0>\leq W\leq\hat W$ 

\n are respectively the Cartan matrix, the root lattice, the root system (with the set of positive roots) and the Weyl group  of the simple Lie algebra $\gothg_0=\gothg^{\chi}$.

\n If $\gothg$ is of type $X_{\tilde n}$ ($X=A,B,C,D,E,F,G,$) $\hat\gothg$ is said to be of type $X_{\tilde n}^{(k)}$.

\n Finally $\hat P=\oplus_{i\in I_0}\Z\lambda_i$ ($<\lambda_i|\alpha_j>=\delta_{ij}\tilde d_i$) 
is the sublattice of $$Hom(Q_0,\Z)\subseteq Hom(Q_0,\Z)\oplus Hom(\Z\delta,\Z)=Hom(Q_0\oplus\Z\delta=Q,\Z)$$ such that $\hat W=\hat P\rtimes W_0$ and $\tilde d_i=\begin{cases}1&{\rm{if}}\ k=1\ {\rm{or}}\ X_{\tilde n}^{(k)}=A_{2n}^{(2)}\cr d_i&{\rm{otherwise}}
\end{cases};$ recall that for all $\tilde\lambda\in \hat P, \alpha\in Q$ $\tilde\lambda(\alpha)=\alpha-<\tilde\lambda|\alpha>\delta$ and denote by $\lambda$ the weight $\lambda=\lambda_1+...+\lambda_n$, by $N$ the length of $\lambda$, by $N_i$ the length of $\lambda_i$ ($i\in I_0$).
\vskip .3 truecm
\rem{\rsy}{}
\n The structure of the set of positive roots with multiplicities $\hat\Phi$ is the following (see \kac):

$\hat\Phi=\Phi_+^{{\rm{re}}}\cup\hat\Phi_+^{{\rm{im}}}$ with

$\Phi_+^{{\rm{re}}}=\{r\delta+\alpha\in Q_+|\alpha\in\Phi_0,r\in\Z\ {\rm{such\ that}}\ \tilde d_{\alpha}|r\}\cup\Phi_2$,

$\hat\Phi_+^{{\rm{im}}}=\{(r\delta,i)|(i,r)\in\iz,r>0
\}$, 

\n where

$\tilde d_{w(\alpha_i)}=\tilde d_i$ for $w\in W_0$, $i\in I_0$,

$\Phi_2=\begin{cases}\{(2r+1)\delta+2\alpha|r\in\N, \alpha\in\Phi_0\ {\rm{such\ that}}\ (\alpha|\alpha)=2\}&{\rm{in\ case\ }}A_{2n}^{(2)}\cr\emptyset&{\rm{otherwise}},\end{cases}$ 

$\iz=\{(i,r)\in I_0\times\Z|\tilde d_i|r\}$.

\vskip .3 truecm

\nota{\gent\ }{(see \beck\ and \damcina)}
\n $\iota:\Z\to I$ and $\Z\ni r\mapsto w_r\in W$ are defined by the following conditions:

i) $w_r=\begin{cases}s_{\iota_1}\cdot s_{\iota_2}\cdot...\cdot s_{\iota_{r-1}}&{\rm{if}}\ r\geq 1\cr s_{\iota_0}\cdot s_{\iota_{-1}}\cdot...\cdot s_{\iota_{r+1}}&{\rm{if}}\ r\leq 0.
\end{cases}$


ii) for all $r=1,...,n$ there exists $\tau_r\in{\cal T}$ such that $$\lambda_1+...+\lambda_r=\lambda_1\cdot...\cdot\lambda_r=s_{\iota_1}\cdot...\cdot s_{\iota_{N_1+...+N_r}}\tau_r\in \hat W;$$ 

iii) $\iota_{N+r}=\tau_n(\iota_r)$ for all $r\in\Z$.

\n The bijection $\Z\ni r\mapsto\beta_r=w_r(\alpha_{\iota_r})\in\Phi_+^{{\rm{re}}}$ induces a total ordering $\preceq$ on $\hat\Phi_+$ defined by
$$\beta_r\preceq\beta_{r+1}\preceq(\tilde m\delta,i)\preceq(m\delta,j)\preceq(m\delta,i)\preceq\beta_{s-1}\preceq\beta_s$$
$$\forall r\geq 1,\ s\leq 0,\ \tilde m>m>0,\ j\leq i\in I_0\ \ ({\rm{choosing\ any\ ordering\ \leq\ of\ }}I_0).$$

(The reverse ordering has the same properties, see \damcina).

\vskip .3 truecm

\nota{\genno}{}
\n i) Consider the ring $\Z[x,x^{-1}]$. Then for all $m,r\in\Z$ the elements $[m]_x$, $[m]_x!$ ($m\geq 0$) and ${m\brack r}_x$ ($m\geq r\geq 0$) are defined respectively by $[m]_x={x^m-x^{-m}\over x-x^{-1}}$, $[m]_x!=\prod_{s=1}^m[s]_x$ and ${m\brack r}_x={[m]_x!\over[r]_x![m-r]_x!}$,
which all lie in $\Z[x,x^{-1}]$. 

\n ii) Consider the field $\C(q)$ and, given $v\in\C(q)\setminus\{0\}$, the natural homomorphism $\Z[x,x^{-1}]\to\C(q)$ determined by the condition $x\mapsto v$; then for all $m,r\in\Z$ the elements $[m]_v$, $[m]_v!$ ($m\geq 0$) and ${m\brack r}_v$ ($m\geq r\geq 0$) denote the images in $\C(q)$ respectively of the elements $[m]_x$, $[m]_x!$ and ${m\brack r}_x$ . 

\n iii) 
For all $i\in I_0$ we denote by $q_i$ the element $q_i=q^{d_i}\in\C(q)$.

\vskip .3 truecm
\nota{\nbr}{}
\n Consider a $\Z[q^{\pm 1}]$-algebra $U$, elements $u,v\in U$ and $r\in\Z$. The $q$-bracket
$[u,v]_{q^r}$ denotes the element $[u,v]_{q^r}=uv-q^rvu$.

\n Remark that the specialization at 1 of $[u,v]_{q^r}$ (the image of $[u,v]_{q^r}$ in the $\Z$-algebra $U/(q-1)U$) is the classical bracket $[u,v]=uv-vu$.

\vskip .5 truecm
\noindent {\bf{\preldj.\ {\bf PRELIMINARIES: $\udj$.}}}
\vskip .5truecm

\n In this section we recall the definition and the structures of the Drinfeld-Jimbo presentation  $\udj
$ of the affine quantum algebras
(see \drr\ and \jm, and also \beck, \damcina, \levsb,\lusz).

\vskip .3 truecm
\ddefi{\QNTALG
}
\noindent The Drinfeld-Jimbo presentation of the affine quantum algebra of type $X_{\tilde n}^{(k)}$ is the $\C(q)$-algebra $\udj=\udj(X_{\tilde n}^{(k)})$ generated by
\vskip .3 truecm
$$\{E_i,F_i,K_i^{\pm 1}|i\in I\}$$ 
with relations: 
$$K_iK_i^{-1}=1=K_i^{-1}K_i,\ \ K_iK_j=K_jK_i\ \ 
\forall i,j\in I,\leqno{(KK)}$$ 
$$K_iE_j=q_i^{a_{ij}}E_jK_i,\ \   K_iF_j=q_i^{-a_{ij}}F_jK_i
\ \ \ 
\forall i,j\in I,\leqno{(KEF)}$$ 
$$[E_i,F_j]=\delta_{ij}{K_i-K_i^{-1}\over q_i-q_i^{-1}}\ \ 
\forall i,j\in I,\leqno{(EF)}$$ 
$$\sum_{u=0}^{1-a_{ij}}(-1)^u{1-a_{ij}\brack u}_{q_i} 
E_i^uE_jE_i^{1-a_{ij}-u}=0\ \ \forall i\neq j\in I,\leqno{(SE)}$$ 
$$\sum_{u=0}^{1-a_{ij}}(-1)^u{1-a_{ij}\brack u}_{q_i} 
F_i^uF_jF_i^{1-a_{ij}-u}\ \ \forall i\neq j\in I.\leqno{(SF)}$$ 

\rem{\qremu}{}
\noindent Recall that  
$\udj$ is endowed with the following structures: 

\noindent i) the $Q$-gradation $\udj=\oplus_{\alpha\in Q}\u_{q,\alpha}^{DJ}$ determined by the
conditions:
$$E_i\in\u_{q,\alpha_i}^{DJ},\ \ F_i\in\u_{q,-\alpha_i}^{DJ},\ \ K_i^{\pm
1}\in\u_{q,0}^{DJ}\ \ \forall i\in I;\ \ \  
\u_{q,\alpha}^{DJ}\u_{q,\beta}^{DJ}\subseteq\u_{q,\alpha+\beta}^{DJ}\ \ \forall\alpha,\beta\in Q;$$

\noindent ii) the triangular decomposition: 
$$\udj\cong\uq^{DJ,-}\otimes\uq^{DJ,0}\otimes\uq^{DJ,+}\cong\uq^{DJ,+}\otimes\uq^{DJ,0}\otimes\uq^{DJ,-},$$ where 
$\uq^{DJ,-}$, $\uq^{DJ,0}$ and $\uq^{DJ,+}$ are the subalgebras of $\udj$ 
generated respectively by 
$\{E_i|i\in I\}$, $\{K_i^{\pm 1}|i\in I\}$ and 
$\{F_i|i\in I\}$; in particular $$\u_{q,\alpha}^{DJ}\cong\bigoplus_{\beta,\gamma\in Q_+:\atop \gamma-\beta=\alpha}\u_{q,-\beta}^{DJ,-}\otimes\uq^{DJ,0}\otimes\u_{q,\gamma}^{DJ,+}\ \  \forall \alpha\in Q$$ 
where $\u_{q,\alpha}^{DJ,\pm}=\u_{q,\alpha}^{DJ}\cap\u_{q}^{DJ,\pm}$ is finite-dimensional $\forall \alpha\in Q$;

\n remark also that if $\tilde\uq^{DJ,-}=\oplus_{\alpha\in Q_+}\tilde\u_{q,-\alpha}^{DJ,-}$ with $\tilde\u_{q,-\alpha}^{DJ,-}=\u_{q,-\alpha}^{DJ,-}K_{\alpha}$, we have that $\tilde\uq^{DJ,-}$ is a graded subalgebra of $\udj$ and the triangular decomposition can be formulated also as 
$$\udj\cong\tilde\uq^{DJ,-}\otimes\uq^{DJ,0}\otimes\uq^{DJ,+}.$$
iii) the $\C$-anti-linear anti-involution 
$\Omega:\udj\rightarrow\udj$ defined by 
$$\Omega(q)= q^{-1},\ \ \Omega(E_i)=F_i,\ \  
\Omega(F_i)=E_i,\ \ 
\Omega(K_i)= K_i^{-1}\ \forall i\in I;$$ 

\n iv) the extended braid group action  
defined by 
$$T_{s_i}(K_j)=K_j
K_i^{-a_{ij}}
\ \ \forall i,j\in I,$$
$$T_{s_i}(E_i)=-F_iK_i,\ \ T_{s_i}(F_i)=-K_i^{-1}E_i \ \ \forall i\in I,$$
$$T_{s_i}(E_j)=\sum_{r=0}^{-a_{ij}}(-1)^{r-a_{ij}}q_i^{-r}
E_i^{(-a_{ij}-r)}E_jE_i^{(r)},\ \ T_{s_i}(F_j)=\Omega(T_{s_i}(E_j)) \ \forall i\neq j\in I$$
where $\forall m\in\N$ $E_i^{(m)}={E_i^m\over[m]_{q_i}!}$,
and 
$$T_{\tau}(K_i)=K_{\tau(i)},\ \ T_{\tau}(E_i)=E_{\tau(i)},\ \ T_{\tau}(F_i)=F_{\tau(i)}\ \ \forall\tau\in Aut(\Gamma),\ i\in I_0;$$

\n v) positive and negative root vectors $E_{\alpha}\in\u_{q,\alpha}^{DJ,+}$ and $F_{\alpha}=\Omega(E_{\alpha})\in\u_{q,-\alpha}^{DJ,-}$ ($\alpha\in\hat{\Phi}_+$) such that 
$E_{\beta_r}=T_{w_r}(E_{\iota_r})$ if $r\geq 1$, $E_{\beta_r}=T_{w_r^{-1}}^{-1}(E_{\iota_r})$ if $r\leq 0$,  and
$$exp\Big((q_i-q_i^{-1})\sum_{r>0} E_{(\tilde d_ir\delta,i)}u^r\Big)=1-(q_i-q_i^{-1})\sum_{r>0} \tilde E_{(\tilde d_ir\delta,i)}u^r$$
where
$\tilde E_{(\tilde d_ir\delta,i)}=-E_{\tilde d_ir\delta-\alpha_i}E_i+q_i^{-2}E_iE_{\tilde d_ir\delta-\alpha_i}$ if $r>0$, $i\in I_0$.

(Remark that in \realdrel\ $E_{(\tilde d_ir\delta,i)}$ was confused with $\tilde E_{(\tilde d_ir\delta,i)}$ by a misprint.)

\vskip .3 truecm
\rem{\grgr}{}
\n We have that: 

\n i) $\Omega\compo T_w=T_w\compo\Omega$ $\forall w\in \hat W$;

\n ii) $T_w(\u_{q,\alpha}^{DJ})=\u_{q,w(\alpha)}^{DJ}$ $\forall w\in\hat{W}$, $\alpha\in Q$; 

\n iii) $T_w(E_i),T_{w^{-1}}^{-1}(E_i)\in\uq^{DJ,+}$ if $w\in \hat W$
and $i\in I$ are such that $w(\alpha_i)\in Q_+$;

\n iv) $T_w(E_i)\in\tilde\uq^{DJ,-}$ if $w\in \hat W$
and $i\in I$ are such that $w(\alpha_i)\in -Q_+$;

\n v) $T_w(E_i)=E_j$ if $w\in \hat W$
and $i\in I$ are such that $w(\alpha_i)=\alpha_j$;

\n vi) $T_{\lambda}(E_{\alpha})=\begin{cases}E_{\lambda(\alpha)}=E_{\alpha-<\lambda|\alpha>\delta}&
{\rm{if}}\ \lambda(\alpha)\in Q_+\cr
-F_{-\lambda(\alpha)}K_{-\lambda(\alpha)}=-F_{<\lambda|\alpha>\delta-\alpha}K_{<\lambda|\alpha>\delta-\alpha}&
{\rm{otherwise}};\end{cases}$


\n vii) $E_{r\tilde d_i\delta+\alpha_i}=T_{\lambda_i}^{-r}(E_i)$ 
$\forall r\in\N, i\in I_0$; 

\n viii) $T_{\lambda_i}(E_{(\tilde d_jr\delta,j)})=E_{(\tilde d_jr\delta,j)}$ for all $i,j\in I_0$, $r>0$;

\n ix) $\{K_{\alpha}=\prod_{i\in I}K_i^{m_i}|\alpha=\sum_{i\in I}m_i\alpha_i\in Q\}$ is a basis of $\uq^{DJ,0}$; 

\n x) $\{E(\gamma)=E_{\gamma_1}\cdot...\cdot E_{\gamma_M}|M\in\N,\ \gamma=(\gamma_1\preceq...\preceq\gamma_{M}),\gamma_{h} \in\hat{\Phi}_+\forall h=1,...,M\}$
is a basis of $\uq^{DJ,+}$ (PBW-basis);

\n xi) $\forall\alpha\prec\beta\in\hat\Phi_+$ $E_{\beta}E_{\alpha}-q^{(\alpha|\beta)}E_{\alpha}E_{\beta}$ is a linear combination of 
$\{E(\gamma)|\gamma=(\gamma_1\preceq...\preceq\gamma_{M})\in\hat{\Phi}_+^M,\ M\in\N,\ \alpha\prec\gamma_1\}$; for real root vectors the claim can be stated in a more precise way: $\forall\beta_r\prec\beta_s\in\hat\Phi_+$ $E_{\beta_s}E_{\beta_r}-q^{(\beta_r|\beta_s)}E_{\beta_r}E_{\beta_s}$ is a linear combination of 
$\{E(\gamma)|\gamma=(\gamma_1\preceq...\preceq\gamma_{M})\in\hat{\Phi}_+^M,\ M\in\N,\ \beta_r\prec\gamma_u\prec\beta_s\ \forall u=1,...,M\}$
(Levendorskii-Soibelman formula).

\vskip .3 truecm
\rem{\qremd}{}
\n The $\a$-subalgebra $\ua$ of $\udj$ generated by $\{E_i,F_i,K_i^{\pm 1}|i\in I\}$ is an integer form of $\udj$:

\n i) $\udj=\C(q)\otimes_{\a}\ua$;

\n ii) $\ua$ is a free $\a$-module; 

\n Moreover:

\n iii) $\ua$ is $T_{s_i}$-stable for all $i\in I$ and $T_{\tau}$-stable for all $\tau\in{\cal T}$: it contains all the root vectors;

\n iv) the subalgebra $\udjp=\ua\cap\uq^{DJ,+}$ is the $\a$-algebra generated by $\{E_i|i\in I\}$ with 
relations $(SE)$; it is a free $\a$-module;

\n v) $\udjp=\oplus_{\alpha\in Q_+}{\cal U}_{\a,\alpha}^{DJ,+}$ where ${\cal U}_{\a,\alpha}^{DJ,+}=\udjp\cap{\cal U}_{q,\alpha}$ is free of finite rank over $\a$;

\n vi) the subalgebra $\u_{\a}^{DJ,0}=\ua\cap\uq^{DJ,0}$ is the commutative $\a$-algebra
$\u_{\a}^{DJ,0}=\a\Big[K_i,{K_i-K_i^{-1}\over q_i-q_i^{-1}}|i\in I\Big]/\Big(K_i\big(K_i-(q_i-q_i^{-1}){K_i-K_i^{-1}\over q_i-q_i^{-1}}\big)|  i\in I\Big)$; it is a free $\a$-module;

\n vii) $\ua\Big/(q-1,K_i-1|i\in I)\cong{\cal U}(\hat{\gothg})$;

\n viii) for all $i\in I$ $T_{s_i}$ induces $\tilde T_{s_i}:{\cal U}(\hat{\gothg})\to{\cal U}(\hat{\gothg})$ and  
$\tilde T_{s_i}\big|_{\hat{\gothg}}\in Aut_{Lie}(\hat{\gothg})$: the image of all the root vectors lies in $\hat{\gothg}$.

\vskip .5 truecm
\noindent {\bf{\preldr.\ {\bf PRELIMINARIES: $\udr$.}}}
\vskip .5truecm

\n The Drinfeld realization $\udr$ of the affine quantum algebras was introduced in \drld, and its defining relations were simplified in \realdrel\ thanks to the ($q$-)commutation with the generators $X_{i,r}^{\pm}$, $H_{i,r}$.

\n Both the original and the simplified sets of relations are useful in this paper: while studying the positive subalgebra $\uq^{Dr,+}$, which contains neither $X_{i,r}^-$ nor $H_{i,r}$, the set of relations given by Drinfeld is the most natural to deal with, and is finally proved to provide a complete set of relations defining $\uq^{Dr,+}$ (see theorem \tttc,i)); viceversa, specializing at 1 the whole $\udr$ provides a presentation of the affine Kac-Moody algebras in terms of the generators $\{x_{i,r}^{\pm},h_{i,r},c\}$, whose relations can be deduced from the simplified relations defining $\udr$ (see theorem \tttt,iv)).

\noindent In this section we recall: the definition of $\udr$ through the simplified relations given in \realdrel\ (definition \drdef); the relations given by Drinfeld (\drld) involving just the positive generators $X_{i,r}^+$'s and holding in $\uq^{Dr,+}$ (notation \drnotu\ and remark \dremd); the structures defined on $\udr
$ ($Q$-gradation, (anti)automorphisms, first remarks about the triangular decomposition).
\vskip .3 truecm
\ddefi{\drdef}
\noindent The Drinfeld realization of the affine quantum algebra of type $X_{\tilde n}^{(k)}$ is the $\C(q)$-algebra $\udr=\udr(X_{\tilde n}^{(k)})$ generated by $$
C^{\pm 1},\ \ \ k_i^{\pm 1}\ (i\in I_0),\ \ \ X_{i,r}^{\pm}\ ((i,r)\in I_0\times\Z)
$$ with relations
$$X_{i,r}^{\pm}=0\ \ \forall (i,r)\in (I_0\times\Z)\setminus I_{\Z},\leqno{(ZX^{\pm})} $$
$$[C,x]=0\ \ \forall x,\ \ \ k_ik_j=k_jk_i\ \ \ (i,j\in I_0),\leqno{(CUK)}$$
$$CC^{-1}=1,\ \ \ k_ik_i^{-1}=1=k_i^{-1}k_i\ \ \ (i\in I_0),\leqno{(CK)}$$
$$k_iX_{j,r}^{\pm}=q_i^{\pm a_{ij}}X_{j,r}^{\pm}k_i\ \ (i\in I_0,\ (j,r)\in I_0\times\Z),\leqno{(KX^{\pm})}$$
$$[X_{i,r}^+,X_{j,s}^-]=\begin{cases} \delta_{i,j}
{C^{-s}k_i\tilde H_{i,r+s}^+ -C^{-r}k_i^{-1}\tilde H_{i,r+s}^-\over q_i-q_i^{-1}}&{\rm {if}}\ \tilde  d_j|s\cr 0&{\rm {otherwise}}
\end{cases} 
\leqno{(XX)}$$
$$ [H_{i,r},X_{j,s}^{\pm}]=\pm b_{ijr}C^{r\mp|r|\over 2}X_{j,r+s}^{\pm}\ \ ((i,r), (j,s)\in\iz,\ \tilde d_i\leq|r|\leq\tilde d_{ij}),\leqno{(HXL^{\pm})}$$
$$[X_{i,r\pm 1}^{\pm},X_{i,r}^{\pm}]_{q^2}=0\ \ (X_{\tilde n}^{(k)}=A_1^{(1)}),\leqno{(X1_{const}^{\pm})}$$
$$[[X_{i,r\pm 1}^{\pm},X_{i,r}^{\pm}]_{q^2},X_{i,r}^{\pm}]_{q^4}=0\ \ (X_{\tilde n}^{(k)}=A_2^{(2)}),\leqno{(X3_{const}^{\pm})}$$
$$\sum_{u=0}^{1-a_{ij}}(-1)^u{1-a_{ij}\brack u}_{q_i}(X_{i,r}^{\pm})^uX_{j,s}^{\pm}(X_{i,r}^{\pm})^{1-a_{ij}-u}=0\ \ \ (n>1),\leqno{(S_{const}^{\pm})}$$
where
$\tilde H_{i,r}^{\pm}$, 
$H_{i,r}$, 
$b_{ijr}$ and $\tilde d_{ij}$ are defined as follows:
$$\tilde H_{i,r}^{\pm}=\begin{cases}1&{\rm{if}}\ r=0\cr
\pm(q_i-q_i^{-1})C^{{r\mp r\over 2}}k_i^{\mp 1}[X_{i,r}^+,X_{i,0}^-]&{\rm{if}}\ \pm r>0\cr0&{\rm{if}}\ \pm r<0;\end{cases}$$
$$\sum_{r\in\Z}\tilde H_{i,\pm r}^{\pm}u^r=exp\left(\pm(q_i-q_i^{-1})
\sum_{r>0}H_{i,\pm r}u^r\right);$$
$$b_{ijr}=\begin{cases}  0&{\rm{if}}\ \tilde d_{i,j}\not|r\cr
{[2r]_q(q^{2r}+(-1)^{r-1}+q^{-2r})\over r}&{\rm{if}}\ (X_{\tilde n}^{(k)},d_i,d_j)=(A_{2n}^{(2)},1,1)\cr
{[\tilde r a_{ij}]_{q_i}\over\tilde r}&{\rm{otherwise,\ with\ }}
\tilde r={r\over \tilde d_{i,j}};\end{cases} $$
$$\tilde d_{ij}=max\{\tilde d_i,\tilde d_j\}.$$

\vskip .3 truecm
\nota{\drnot}{}

\noindent In $\udr$:

\noindent i) $\uq^{Dr,0}$ denotes the $\C(q)$-subalgebra generated by $\{C^{\pm 1},k_i^{\pm 1},H_{i,r}|i\in I_0,r\neq 0\}$, or, equivalently, the $\C(q)$-subalgebra generated by $\{C^{\pm 1},k_i^{\pm 1},\tilde H_{i,r}^{\pm}|i\in I_0,r\in\Z\}$;

\noindent ii) $\uq^{Dr,0,0}$, $\uq^{Dr,0,+}$ and $\uq^{Dr,0,-}$ denote the $\C(q)$-subalgebras generated respectively by $\{C^{\pm 1}, k_i^{\pm 1}|i\in I_0\}$, by $\{H_{i,r}|i\in I_0,r>0\}$ (or by $\{\tilde H_{i,r}^+|i\in I_0,r\in\Z\}$) and by $\{H_{i,r}|i\in I_0,r<0\}$ (or by $\{\tilde H_{i,r}^-|i\in I_0,r\in\Z\}$);

\noindent iii) $\uq^{Dr,+}$ and 
$\uq^{Dr,-}$ denote the $\C(q)$-subalgebras generated respectively by $\{X_{i,r}^+|i\in I_0,r\in\Z\}$ and 
by $\{X_{i,r}^-|i\in I_0,r\in\Z\}$;

\noindent iv) $\uq^{Dr,+,+}$ and $\uq^{Dr,-,-}$ denote the $\C(q)$-subalgebras generated respectively by $\{X_{i,r}^+|i\in I_0,r\geq\Z\}$ and by $\{X_{i,r}^-|i\in I_0,r\leq\Z\}$;

\noindent v) given $\alpha\in Q$, $\u_{q,\alpha}^{Dr}$ denotes the $\alpha$-homogeneous component of $\udr$: 
$\udr=\oplus_{\alpha\in Q}\u_{q,\alpha}^{Dr}$ where $C^{\pm 1}, k_i^{\pm 1}\in\u_{q,0}$, $X_{i,r}^{\pm}\in\u_{q,r\delta\pm\alpha_i}^{Dr}$; $\u_{q,\alpha}^{Dr,*}=\uq^{Dr,*}\cap\u_{q,\alpha}^{Dr}$.

\vskip .3 truecm
\rem{\drrem}{}

\noindent i) $\uq^{Dr,0,0}\subseteq\u_{q,0}^{Dr}$;

\noindent ii) $\uq^{Dr,0}\subseteq\oplus_{m\in\Z}\u_{q,m\delta}^{Dr}$;

\noindent iii) $\uq^{Dr,+}\subseteq\C(q)\oplus\big(\oplus_{m\in\Z,\alpha\in Q_{0,+}\setminus\{0\}}\u_{q,m\delta+\alpha}^{Dr}\big)$;

\noindent iv) $\uq^{Dr,+,+}\subseteq\C(q)\oplus\big(\oplus_{m\in\N,\alpha\in Q_{0,+}\setminus\{0\}}\u_{q,m\delta+\alpha}^{Dr}\big)$;

\n v) for all $\alpha\in Q_{0,+}$, $m\in\Z$ ${\cal U}_{q,m\delta+\alpha}^{Dr,+,+}$ is finite-dimensional, while ${\cal U}_{q,m\delta+\alpha}^{Dr,+}$ is in general not finite-dimensional.

\vskip .3 truecm
\ddefi{\drom}

\noindent i) $\Omega:\udr\to\udr$ is the $\C$-anti-linear anti-involution defined by 
$$q\mapsto q^{-1},\ \ C^{\pm 1}\mapsto C^{\mp 1},\ \ k_i^{\pm 1}\mapsto k_i^{\mp 1},\ \ X_{i,r}^{\pm}\mapsto X_{i,-r}^{\mp},$$
$$\tilde H_{i,r}^{\pm}\mapsto \tilde H_{i,-r}^{\mp},\ \ H_{i,r}\mapsto H_{i,-r}.$$
ii) $t_i:\udr\to\udr$ ($i\in I_0$) is the $\C(q)$-automorphism defined by
$$C^{\pm 1}\mapsto C^{\pm 1},\ \ \ k_j^{\pm 1}\mapsto (k_jC^{-\delta_{ij}\tilde d_i})^{\pm 1},\ \ \ X_{j,r}^{\pm}\mapsto X_{j,r\mp\delta_{ij}\tilde d_i}^{\pm},$$
$$\tilde H_{j,r}^{\pm}\mapsto \tilde H_{j,r}^{\pm}\ \ \ H_{j,r}\mapsto H_{j,r}.$$
iii) ${\cal E}_c:\udr\to\udr$ ($c:I_0\to\{\pm 1\}$
) is the $\C(q)$-automorphism defined by
$${\cal E}_c\big|_{\uq^{Dr,0}}=id_{\uq^{Dr,0}},\ \ \ X_{i,r}^{\pm}\mapsto c_i
X_{i,r}^{\pm}.$$
\vskip .3 truecm
\rem{\drromrem}{}

\noindent i) For all $i,j\in I_0$ we have $\Omega\compo t_i=t_i\compo\Omega$ and $t_i\compo t_j=t_j\compo t_i$;

\noindent ii) $\Omega(\uq^{Dr,0})=\uq^{Dr,0}$, $\Omega(\uq^{Dr,0,0})=\uq^{Dr,0,0}$, $\Omega(\uq^{Dr,0,\pm})=\uq^{Dr,0,\mp}$, $\Omega(\uq^{Dr,\pm})=\uq^{Dr,\mp}$, $\Omega(\u_{q,\alpha}^{Dr})=\u_{q,-\alpha}^{Dr}$;

\n iii) $t_i(\u_{q,\alpha}^{Dr})=\u_{q,\lambda_i(\alpha)}^{Dr}$;

\noindent iv) $t_i(\uq^{Dr,*})=\uq^{Dr,*}$, $t_i(\uq^{Dr,0,*})=\uq^{Dr,0,*}$;

\noindent v) more precisely $t_i\big|_{\uq^{Dr,0,\pm}}=id_{\uq^{Dr,0,\pm}}$ ($*\in\{0,+,-\}$);

\n vi) $t_i^{-1}(\uq^{Dr,+,+})\subseteq\uq^{Dr,+,+}$; 

\n vii) $\uq^{Dr,+}=\cup_{N\in\N}(t_1\cdot...\cdot t_n)^N(\uq^{Dr,+,+})$; 

\n viii) for all $c,\tilde c:I_0\to\{\pm 1\}
$ and for all $i\in I_0$ we have ${\cal E}_c\compo{\cal E}_{\tilde c}={\cal E}_{c\tilde c}$, ${\cal E}_c\compo\Omega=\Omega\compo{\cal E}_c$ and ${\cal E}_c\compo t_i=t_i\compo{\cal E}_c$;

\n ix) ${\cal E}_c(x)=\pm x$ for all $x\in\u_{q,\alpha}^{Dr}$, $\alpha\in Q$.

\vskip .3 truecm
\rem{\uqz}{}

\noindent In $\udr$ we have also (see \drld\ and \realdrel): 
$$H_{i,r}=0\ \ \forall (i,r)\in (I_0\times\Z)\setminus I_{\Z},\leqno{(ZH)} $$
$$[k_i,H_{j,s}]=0\ \ (i\in I_0,\ (j,s)\in I_0\times(\Z\setminus\{0\})),\leqno{(KH)}$$
$$ [H_{i,r},H_{j,s}]=\delta_{r+s,0}b_{ijr}{C^r-C^{-r}\over q_j-q_j^{-1}}\ \ ((i,r),(j,s)\in I_0\times(\Z\setminus\{0\})),\leqno{(HH)}$$
so that:

\noindent i) $\uq^{Dr,0,0}$ is central in $\uq^{Dr,0}$; 

\noindent ii) $\uq^{Dr,0,0}$ is a quotient of $\C(q)[C^{\pm 1},k_i^{\pm 1}|i\in I_0]$ and $\uq^{Dr,0,+}$ is a quotient of $\C(q)[H_{i,r}|i\in I_0,\tilde d_i|r>0]$; 

\noindent iii) the natural homomorphism of $\C(q)$-vector spaces
$$\uq^{Dr,0,-}\otimes_{\C(q)}\uq^{Dr,0,0}\otimes_{\C(q)}\uq^{Dr,0,+}\to\uq^{Dr,0}$$
is surjective.

\vskip .3 truecm
\rem{\dremrelu}{}

\noindent In  $\udr$ we have also (see \drld\ and \realdrel): 
$$ [H_{i,r},X_{j,s}^{\pm}]=\pm b_{ijr}C^{r\mp|r|\over 2}X_{j,r+s}^{\pm}\ \ ((i,r), (j,s)\in I_0\times\Z,
r\neq 0),\leqno{(HX^{\pm})}$$
which, together with the relations $(CUK)$, $(CK)$, $(KX^{\pm})$, $(XX)$, implies that the natural map 
$\uq^{Dr,-}\otimes\uq^{Dr,0}\otimes\uq^{Dr,+}\to\udr$ is surjective.

\rem{\stagg} {}
\n Notice that setting $\tilde\uq^{Dr,-}=\oplus_{\alpha\in Q}K_{\alpha}\u_{q,\alpha}^{Dr,-}(=\oplus_{\alpha\in Q}\u_{q,\alpha}^{Dr,-}K_{\alpha})$ we get that
$$\tilde\uq^{Dr,-}\otimes\uq^{Dr,0}\otimes\uq^{Dr,+}\cong\uq^{Dr,-}\otimes\uq^{Dr,0}\otimes\uq^{Dr,+}.$$

\vskip .3 truecm
\nota{\drnotu}{}

\noindent i) Denote by
$(DR)$ the following relations:

$$[X_{i,r+\tilde d_{ij}}^{+},X_{j,s}^{+}]_{q_i^{a_{ij}}}
+[X_{j,s+\tilde d_{ij}}^{+},X_{i,r}^{+}]_{q_j^{a_{ji}}}=0\ ((i,r),(j,s)\in \iz,\ 
a_{ij}\!\!<\!\!0),\leqno{(XD^{})}$$
$$
\sum_{\sigma\in\sy_2}\sigma.[X_{i,r_1+\tilde d_i}^{+},X_{i,r_2}^{+}]_{q_i^{2}}=0\ \ ((r_1,r_2)\in\Z^2,\ (X_{\tilde n}^{(k)},d_i)\neq(A_{2n}^{(2)},1)),\leqno{(X1^{})}$$
$$
\sum_{\sigma\in\sy_2}\sigma.([X_{i,r_1+ 2}^{+},X_{i,r_2}^{+}]_{q^{2}}-q^{4}[X_{i,r_1+ 1}^{+},X_{i,r_2+ 1}^{+}]_{q^{-6}})=0\leqno{(X2^{})}$$
$\ \ \ \ ((r_1,r_2)\in\Z^2,\ (X_{\tilde n}^{(k)},d_i)=(A_{2n}^{(2)},1)),$
$$\sum_{\sigma\in\sy_3}\sigma.[[X_{i,r_1+ 1}^{+},X_{i,r_2}^{+}]_{q^{ 2}},X_{i,r_3}^{+}]_{q^{ 4}}=0\leqno{(X3)}$$
$\ \ \ \ ((r_1,r_2,r_3)\in\Z^3,\ (X_{\tilde n}^{(k)},d_i)=(A_{2n}^{(2)},1)),$
$$\sum_{\sigma\in\sy_{1-a_{ij}}}\sigma.
[[...[[X_{j,s}^+,X_{i,r_1}^{+}]_{q_i^{-a_{ij}}},X_{i,r_2}^+]_{q_i^{-a_{ij}-2}},...]_{q_i^{a_{ij}+2}}
,X_{i,r_{1-a_{ij}}}^{+}]_{q_i^{a_{ij}}}
=0\leqno{(SUL^{}
)}$$
$$\ \ \ \ (i\neq j\in\! I_0,\ a_{ij}\in\{0,-1\}\ {\rm {if}}\ k\neq 1,
\ r\in
\Z^{1-a_{ij}},\ s\in\Z),$$
$$\sum_{\sigma\in\sy_2}\sigma.[[X_{j,s}^{+},X_{i,r_1+ 1}^{+}]_{q^2},X_{i,r_2}^{+}]=0,\leqno(T2^{})$$
$$\ \ \ \ (i,j\in I_0,\ a_{ij}=-2,\ k=2,\ X_{\tilde n}^{(k)}\neq A_{2n}^{(2)},\ (r_1,r_2)\in\Z^{2},\ s\in\Z),$$
$$\sum_{\sigma\in\sy_2}\sigma.\big((q^2+q^{-2})[[X_{j,s}^{+},X_{i,r_1+ 1}^{+}]_{q^2},X_{i,r_2}^{+}]+\leqno(S2^{})$$
$$+q^2[[X_{i,r_1+ 1}^{+},X_{i,r_2}^{+}]_{q^{ 2}},X_{j,s}^{+}]_{q^{- 4}}\big)=0$$
$$\ \ \ \ (i,j\in I_0,\ a_{ij}=-2,\ X_{\tilde n}^{(k)}= A_{2n}^{(2)},\ (r_1,r_2)\in\Z^{2},\ s\in\Z),$$
$$\sum_{\sigma\in\sy_2}\sigma.((q^{2}+1)[[X_{j,s}^{+},X_{i,r_{1}+ 2}^{+}]_{q^{3}},X_{i,r_{2}}^{+}]_{q^{-1}}+\leqno{(T3^{})}$$
$$+[[X_{j,s}^{+},X_{i,r_{1}+ 1}^{+}]_{q^{3}},X_{i,r_{2}+ 1}^{+}]_{q})=0$$
$$\ \ \ \ (i,j\in I_0,\ a_{ij}=-3,\ k=3,\ (r_1,r_2)\in\Z^{2},\ s\in\Z).$$
\n ii) Denote by $(S)$ the relations
$$\sum_{\sigma\in\sy_{1-a_{ij}}}\sigma.
[...[[X_{j,s}^+,X_{i,r_1}^{+}]_{q_i^{-a_{ij}}},X_{i,r_2}^+]_{q_i^{-a_{ij}-2}},...
,X_{i,r_{1-a_{ij}}}^{+}]_{q_i^{a_{ij}}}
=0\leqno{(S
)}$$
$$\ \ \ \ (i\neq j\in\! I_0,
\ r\in
\Z^{1-a_{ij}},\ s\in\Z),$$

\noindent iii) Denote by $(U3)$ the relations
$$\sum_{\sigma\in\sy_{3}}\sigma.
[[[X_{j,s}^+,X_{i,r_1+1}^{+}]_{q^{3}},X_{i,r_2}^+]_{q^{}},X_{i,r_{3}}^{+}]_{q^{-1}}
=0\leqno{(U3)}$$
$$\ \ \ \ (i,j\in I_0,\ a_{ij}=-3,\ k=3,\ (r_1,r_2)\in\Z^{2},\ s\in\Z).$$

\noindent iv) Denote by $(ZX_+^+)$, respectively $(DR_+)$,  $(S_+)$ and $(U3_+)$, the relations of $(ZX^+)$, respectively $(DR)$, $(S)$ and $(U3)$, involving just elements $X_{i,r}^+$ with $r\geq 0$ (see definition \drdef).

\vskip .3 truecm
\rem{\dremd}{}

\n i) The relations $(DR)$, $(S)$ and $(U3)$ hold in $\udr$ (see \drld\ and \realdrel); 

\n ii) the relations $(SUL)$ and $(SUL_+)$ depend respectively on the relations $(S)$ and $(S_+)$;

\n iii) the relations $(S)$ and $(U3)$ depend on the relations $(DR)$ (see \realdrel);

\n iv) in the algebra generated by $\{X_{i,r}^+|i\in I_0, r\in\N\}$ the relations $(S_+)$ and $(U3_+)$ do not depend on the relations $(DR_+)$ (it is enough to compare the degrees of the relations $(S_+)$ and $(U3_+)$ with those of the relations $(DR_+)$ remarking that the algebra generated by $\{X_{i,r}^+|i\in I_0, r\in\N\}$ is $Q_{0,+}\oplus\N\delta$-graded).

\vskip .5 truecm
\noindent {\bf{\prelpsi\ {\bf PRELIMINARIES: $\psi$.}}}
\vskip .5truecm
\noindent In this section we recall the homomorphism $\psi:\udr\to\udj$ and some of its properties (see \beck\ and \damcina).

\vskip .3truecm

\ddefi{\defpsi
}
\n $\psi=\psi_{X_{\tilde n}^{(k)}}:{\cal U}^{Dr}_q(X_{\tilde n}^{(k)})\to\udj(X_{\tilde n}^{(k)})$ is the $\C(q)$-algebra homomorphism 
defined on the generators as follows: 
$$C^{\pm1}\mapsto K_{\delta}^{\pm1},\ \ \ 
k_i^{\pm1}\mapsto K_i^{\pm1}\ \ (i\in I_0),$$
$$X_{i,\tilde d_i r}^+\mapsto o(i)^rT_{\lambda_i}^{-r}(E_i),\ \ \ 
X_{i,\tilde d_i r}^-\mapsto o(i)^rT_{\lambda_i}^r(F_i)\ \ (i\in I_0, r\in\Z),$$
$$H_{i,\tilde d_i r}\mapsto \begin{cases}  o(i)^rE_{(\tilde d_i r\d,i)}&{\rm{if\ }} r>0\cr
o(i)^rF_{(-\tilde d_i r\d,i)}&{\rm{if\ }} r<0\end{cases} \ \ (i\in I_0, r\in\Z\setminus\{0\}),$$
where $o:I_0\to\{\pm1\}$ is a map such that: 

\n a) $a_{ij}\neq 0\Rightarrow o(i)o(j)=-1$ (see \beck\ for the untwisted case);

\n b) in the twisted case different from $A_{2n}^{(2)}$ $ a_{ij}=-2\Rightarrow o(i)=1$ (see \realdrel).

\vskip .3 truecm

\rem{\repsi}{\ (see \realdrel)}
\n i) $\psi$ preserves the gradation, that is $\psi(\u_{q,\alpha}^{Dr})=\u_{q,\alpha}^{DJ}$ for all $\alpha\in Q$;

\n ii) $\psi\compo\Omega=\Omega\compo\psi$;

\n iii) $\psi\compo {\cal E}_{o_i}\compo t_i=T_{\lambda_i}\compo\psi$ $\forall i\in I_0$, where $o_i(j)=\begin{cases}o(i)&{\rm{if}}\ j=i\cr 1&{\rm{otherwise}};\end{cases}$

\n iv) $\psi$ is surjective.

\vskip .3 truecm

\prop{\copsi}
\n Let us compare $\udr$ and $\udj$ through $\psi$; then the PBW basis of $\udj$ and remark \uqz, ii) and iii) imply that:

\n i) $\uq^{Dr,0,0}\cong\C(q)[C^{\pm 1},k_i^{\pm 1}|i\in I_0]$ and $\psi\big|_{\uq^{Dr,0,0}}:\uq^{Dr,0,0}\to\uq^{DJ,0}$ is an isomorphism; 

\n ii) $\uq^{Dr,0,+}\!\!\cong\!\!\C(q)[H_{i,r}|i\in I_0, \tilde d_i|r\!>\!0]$ and $\uq^{Dr,0,-}\!\!\cong\!\!\C(q)[H_{i,r}|i\in I_0, \tilde d_i|r\!<\!0]$; 

\n iii) the composition
$$\uq^{Dr,0,-}\otimes_{\C(q)}\uq^{Dr,0,0}\otimes_{\C(q)}\uq^{Dr,0,+}\to\uq^{Dr,0}\hookrightarrow\udr\overset{\psi}{\to}\udj$$
is injective; 

\n iv) $\uq^{Dr,0,-}\otimes_{\C(q)}\uq^{Dr,0,0}\otimes_{\C(q)}\uq^{Dr,0,+}\cong\uq^{Dr,0}$;

\n v) $\psi\big|_{\uq^{Dr,0}}:\uq^{Dr,0}\to\udj$ is injective.
\vskip .3 truecm

\vskip .5 truecm
\noindent {\bf{\trngfindim.\ {\bf REDUCTION to a FINITE DIMENSIONAL SITUATION and TRIANGULAR DECOMPOSITION.}}}
\vskip .5truecm

\n The aim of this paper is proving that $\psi$ is an isomorphism, i.e. that it~is injective (since it is surjective).
The strategy is reducing to studying the restriction of $\psi$ to finitely generated $\a$-submodules of $\udr$, so that the specialization argument described in the introduction (proposition \specgen) can be applied.

\n The first step in this direction would be restricting to the $Q$-homogeneous components ${\cal U}_{q,\alpha}^{Dr}$, which are though far from being finite-dimensional; in similar situations, for example while studying the Drinfeld-Jimbo presentation of quantum algebras, the triangular decomposition solves this difficulty, because it provides the lower bound $0\in Q$ for the weight of the elements to be considered. 

\n In the Drinfeld realization this simplification is important but not enough: indeed ${\cal U}_{q,\alpha}^{Dr,+}$ is in general not finite-dimensional (see remark \drrem,v)). The same remark suggests to analyze in fact $\uq^{Dr,+,+}$ since it is the direct sum of (its homogeneous) finite-dimensional components.

\n This section is devoted to show that the injectivity of $\psi\big|_{\uq^{Dr,+,+}}$ implies the injectivity of $\psi$.

\n As outlined above, the reduction to this finite-dimensional situation requires the analysis and understanding of the triangular decomposition of $\udr$.

\n By triangular decomposition of $\udr$ we mean the following claim: 
$$\uq^{Dr,-}\otimes_{\C(q)}\uq^{Dr,0}\otimes_{\C(q)}\uq^{Dr,+}\cong\udr.$$
In \hern\ the author proved the triangular decomposition for the quantum affinizations of all symmetrizable quantum algebras: this class of algebras includes the untwisted affine quantum algebras, but does not include the twisted ones. 

\n Here we develop some remarks which show that the injectivity of $\psi\big|_{\uq^{Dr,+,+}}$ implies both the triangular decomposition of $\udr$ and the injectivity of $\psi$.

\n We already noticed that the product 
$\uq^{Dr,-}\otimes_{\C(q)}\uq^{Dr,0}\otimes_{\C(q)}\uq^{Dr,+}\to\udr$ is surjective (see remark \dremrelu): therefore the triangular decomposition is equivalent to the injectivity of this map.

\vskip .3 truecm
\prop{\prpsid}
\n The product map $\psi(\uq^{Dr,-,-})\otimes_{\C(q)}\psi(\uq^{Dr,0})\otimes_{\C(q)}\psi(\uq^{Dr,+,+})\to\udj$
is injective.
\dim
$\psi(\uq^{Dr,+,+})$ is the subalgebra of $\udj$ generated by the root vectors $E_{r\delta+\alpha_i}$ ($i\in I_0$, $r\in\N$), hence, by the Levendorskii-Soibelman formula and the PBW-basis (see remark \grgr, x) and xi)), it is a subspace of the linear span of the ordered monomials in the root vectors $E_{\beta_r}$ with $r\leq 0$. Of course $\psi(\uq^{Dr,-,-})=\Omega(\psi(\uq^{Dr,+,+}))$, hence it is a subspace of the linear span of the ordered monomials in the root vectors $F_{\beta_r}$ with $r\leq 0$. 

\n Recall that
 $\psi(\uq^{Dr,0})\cong\psi(\uq^{Dr,0,-})\otimes\psi(\uq^{Dr,0,0})\otimes\psi(\uq^{Dr,0,+})$ (see proposition \copsi, iv) and v)), and that $\psi(\uq^{Dr,0,+})$ is the subalgebra of $\udj$ generated by the root vectors $E_{(r\delta,i)}$ ($i\in I_0$, $r>0$).
 
 \n Then the triangular decomposition of $\udj$ (see remark \qremu,ii)) and the structure of its PBW-basis (see remark \grgr, x)) imply the assertion, thanks to proposition \copsi, i) and ii).

\vskip .3 truecm
\prop{\prpsit}
\n If $\psi\big|_{\uq^{Dr,+,+}}$ is injective then:

\n i) $\psi$ is injective (that is $\udr\cong\udj$, see remark \repsi,iv)); 

\n ii) $\uq^{Dr,-}\otimes_{\C(q)}\uq^{Dr,0}\otimes_{\C(q)}\uq^{Dr,+}\cong\udr$.
\dim
It is enough to consider the following commutative diagram for all $N\in\N$ (see propositions
\copsi\ and \prpsid\ and remarks \drromrem, vii) and viii), \dremrelu\ and \repsi,iii)): 
$$\xymatrix{&\uq^{Dr,-,-}\otimes\uq^{Dr,0}\otimes\uq^{Dr,+,+}\ar[rr]^{({\cal E}_o\compo t_1\compo...\compo t_n)^N}\ar@{^{(}->}[d]
&{}
&\uq^{Dr,-}\otimes\uq^{Dr,0}\otimes\uq^{Dr,+} \ar@{->>}[d]^{\mu^{Dr}}\\
&\psi(\uq^{Dr,-,-})\otimes\psi(\uq^{Dr,0})\otimes\psi(\uq^{Dr,+,+})\ar@{^{(}->}[d]
&{}
&\udr\ar[d]^{\psi}\\
&\udj\ar[rr]^{T_{\lambda}^N}
&{}
&\udj}$$
\vskip .3 truecm
\rem{\repsit}{}
\n $\psi(\uq^{Dr,+,+})\subseteq\uq^{DJ,+}$.

\n On the other hand $\psi(\uq^{Dr,+})\not\subseteq\uq^{DJ,+}$. More precisely for all $i\in I_0$, $r>0$ 
$\psi(X_{i,-r}^+)\in\tilde\uq^{DJ,-}$ and $\psi(\uq^{Dr,+})\cap\tilde{\cal U}_{q,-r\delta+\alpha_i}^{DJ,-}\neq \{0\}$ if $\tilde d_i|r$.

\n In particular the Drinfeld triangular decomposition that we aim to prove would not correspond to the Drinfeld and Jimbo triangular decomposition, but would give rise to a substantially different decomposition (Drinfeld triangular decomposition). For a comparison between the two decompositions see proposition \tttq. 
\vskip .3 truecm

\lem{\epp}
\n Let $\alpha\in Q_{0,+}$, $r\geq 0$, $i\in I_0$ be such that $r\delta+\alpha\in\Phi^{{\rm{re}}}$, or $(r\delta,i)\in\hat\Phi^{{\rm{im}}}$. Then

\n i) $E_{r\delta+\alpha}\in\psi({\uq^{Dr,+}})$, and if $r>0$ $F_{r\delta-\alpha}K_{r\delta-\alpha}\in\psi({\uq^{Dr,+}})$;

\n ii) $K_{\alpha-r\delta}E_{r\delta-\alpha}\in\psi(\uq^{Dr,-})$ if $r>0$; 

\n iii) $E_{(r\delta,i)}\in\psi(\uq^{Dr,0})$ if $r>0$;

\dim
Let $U\subseteq\udj$ be defined by
$$U=\{x\in\udj|\forall N>>0\ T_{\lambda}^{-N}(x)\in\uq^{DJ,+},\ T_{\lambda}^{N}(x)\in\tilde\uq^{DJ,-}\}.$$
Then:

\n a) $U$ is a $T_{\lambda}^{\pm 1}$-stable $\C(q)$-subalgebra of $\udj$ (obvious).

\n b) $\psi(\uq^{Dr,+,+})\subseteq U$ thanks to remarks \grgr, vi), \drromrem, vi) and \repsit.

\n c) $\psi(\uq^{Dr,+})\subseteq U$ thanks to a) , b) and remark \drromrem, vii).

\n d) $U\!=\!\psi(\uq^{Dr,+})$: consider the identifications induced by the product $$\udj\cong\psi(\uq^{Dr,-})\otimes\psi(\uq^{Dr,0})\otimes\psi(\uq^{Dr,+})\cong$$
$$\cong\psi(\uq^{Dr,-})\otimes\psi(\uq^{Dr,0,-})\otimes\psi(\uq^{Dr,0,0})\otimes\psi(\uq^{Dr,0,+})\otimes\psi(\uq^{Dr,+})$$
and remark that through these isomorphisms $\forall u\in \udj$ $\exists\tilde N\in\Z$ such that for all $N>\tilde N$ $$T_{\lambda}^{-N}(u)\in\psi(\uq^{Dr,-,-})\otimes\psi(\uq^{Dr,0,-})\otimes\psi(\uq^{Dr,0,0})\otimes\psi(\uq^{Dr,0,+})\otimes\psi(\uq^{Dr,+,+});$$
moreover $$\psi(\uq^{Dr,-,-})\otimes\psi(\uq^{Dr,0,-})\subseteq\uq^{DJ,-},$$  
$$\psi(\uq^{Dr,0,0})\subseteq\uq^{DJ,0},$$  
$$\psi(\uq^{Dr,0,+})\otimes\psi(\uq^{Dr,+,+})\subseteq\uq^{DJ,+};$$
hence if $u\in U$ the condition $T_{\lambda}^{-N}(u)\in\uq^{DJ,+}$ for all $N>>0$ and the triangular decomposition of $\udj$ imply that $u\in\psi(\uq^{Dr,0,+})\otimes\psi(\uq^{Dr,+})$;
but then $\forall N>>0$ $$T_{\lambda}^N(u)\in\psi(\uq^{Dr,0,+})\otimes\tilde\uq^{DJ,-},$$ and again since $\psi(\uq^{Dr,0,+})\subseteq\uq^{DJ,+}$, the condition 
$T_{\lambda}^{N}(u)\in\tilde\uq^{DJ,-}$ for all $N>>0$ and the triangular decomposition of $\udj$ imply that $u\in\psi(\uq^{Dr,+})$, which implies the claim.

\n e) $E_{r\delta+\alpha}\in U$ thanks to remark \grgr, vi).

\n f) $F_{r\delta-\alpha}K_{r\delta-\alpha}\in U$ thanks to a) and e), since $F_{r\delta-\alpha}K_{r\delta-\alpha}$ is $T_{\lambda}$-conjugate to any $E_{s\delta+\alpha}$ with $s\geq 0$ such that $<\lambda|\alpha>|r+s$ (see remark \grgr, vi)).

d), e) and f) imply i).

\n Applying $\Omega$ to f) we get ii), while iii) is a straightforward consequence of the definitions.

\vskip .3 truecm
\cor{\pip}
\n $\uq^{DJ,+}\cap\psi(\uq^{Dr,+})$ is the $\C(q)$-linear span of the ordered monomials in the $E_{r\delta+\alpha}$'s with $r\geq 0$, $\alpha\in Q_{0,+}$ such that $r\delta+\alpha\in\Phi^{{\rm{re}}}$.
\dim
Let $U_+$ be the $\C(q)$-linear span of the ordered monomials in the $E_{\beta_r}$'s with $r\leq 0$, $U_{-}$ be the $\C(q)$-linear span of the ordered monomials in the $E_{\beta_r}$'s with $r\geq 1$ and $U_{0}$ be the $\C(q)$-linear span of the monomials in the positive imaginary root vectors.
Then the $PBW$-basis of $\uq^{DJ,+}$ says that $\uq^{DJ,+}\cong U_{-}\otimes U_{0}\otimes U_+$.
But 
$$U_{-}\otimes U_0\subseteq\psi(\uq^{Dr,-})\otimes\psi(\uq^{Dr,0}),\ \ U_+\subseteq\psi(\uq^{Dr,+})$$
and
$$\psi(\uq^{Dr,-})\otimes\psi(\uq^{Dr,0})\otimes\psi(\uq^{Dr,+})\cong\udj,$$
so that 
$\uq^{DJ,+}\cap\psi(\uq^{Dr,+})\subseteq U_+$,
which is the assertion, thanks to lemma \epp, i).

\vskip .5 truecm

\noindent {\bf{\ntgf.\ {\bf INTEGER FORM.}}}
\vskip .5truecm

\n We are reduced to prove that $\psi\big|_{\uq^{Dr,+,+}}:\uq^{Dr,+,+}\to\uq^{DJ,+}$ is injective, and we want to show it through specialization at 1. This requires to pass to integer forms of  $\uq^{Dr,+,+}$ and $\uq^{DJ,+}$ and to their presentations by generators and relations.

\n To this aim we start with some notations, underlining that in this section we work with the ring $\a=\C[q]_{(q-1)}$ (the localization 
of $\C[q]$ at $(q-1)$).

\vskip .3 truecm

\nota{\spnot}{}

\noindent i) $\f_+$ is the $\a$-algebra freely generated by $\{X_{i,r}^+|i\in I_0,r\geq 0\}$;

\noindent ii) 
$\i_+$ is the ideal of $\f_+$ defined by the relations $(ZX_+^+)$, $(DR_+)$, $(S_+)$ and $(U3_+)$ (see notation \drnotu);

\noindent iii) 
$t_+^{{\prime}}:\f_+\to\f_+$ is the $\a$-endomorphism defined by $X_{i,r}^+\mapsto o(i)X_{i,r+\tilde d_i}^+$ (see definitions \drom, ii) and iii) and \defpsi); we also denote 
by $\bar t_+^{{\prime}}$  the $\a$-endomorphism induced by $t_+^{{\prime}}$ on $\f_+/\i_+$.

\vskip .3 truecm
\rem{\firem}{}
\n i) $\f_+$, $\i_+$ and consequently also $\f_+/\i_+$ are all   $Q$-graded 
where the degree of $X_{i,r}^+$ is $\alpha_i+r\delta$;

\n ii) the $\a$-modules $(\f_+)_{\alpha}$ and $(\f_+/\i_+)_{\alpha}$ ($\alpha\in Q$) are finitely generated: they are generated over $\a$ by $$\{X_{i_1,m_1}^+\cdot...\cdot X_{i_h,m_h}^+|i_r\in I_0,m_r\geq 0\ \forall r=1,...,h,\ \sum_{r=1}^h m_r\delta+\alpha_{i_r}=\alpha\};$$

\n iii) the natural map $f_+:\f_+/\i_+
\to\udr$ is well defined (see definition \drdef\ and remark \dremd,i)); 

\n iv) $f_+\compo\bar t_+^{{\prime}}={\cal E}_o\compo t_1^{-1}\compo...\compo t_n^{-1}
\compo f_+$.

\vskip .3 truecm
\rem{\insint}{}
\n Of course $\C(q)\otimes_{\a}f_+(\f_+/\i_+)\overset{\cong}{\to}\uq^{Dr,+,+}$ and $f_+(\f_+/\i_+)$ is an integer form of $\uq^{Dr,+,+}$: indeed $f_+(\f_+/\i_+)$ is direct sum of finitely generated $\a$-submodules of a $\C(q)$-vector space, hence it is free over $\a$.

\n In particular a $\C(q)$-linear map defined on $\uq^{Dr,+,+}$ is injective if and only if its restriction to $f_+(\f_+/\i_+)$ is  injective.
\vskip .3 truecm

\cor{\pfp}
\n If $\psi\compo f_+$ is injective then:

\n i) $f_+$ is injective, hence $\f_+/\i_+$ is an integer form of $\uq^{Dr,+,+}$, see remark \insint;

\n ii) $\psi\big|_{f_+(\f_+/\i_+)}$ is injective (then so are $\psi\big|_{\uq^{Dr,+,+}}$ and $\psi$, see proposition \prpsit, i) and remark \insint).

\vskip .3 truecm
\rem{\imredu}{}
\n The image of $\psi\compo f_+$ is contained in $\udjp$. 

\n Indeed $\psi(X_{i,r}^+)\in\ua\cap\uq^{DJ,+}=\udjp$ if $r\geq 0$ (see definition \defpsi\ and remarks \grgr, iii) and \qremd, iii)).
\vskip.3 truecm
\nota{\uffnot}{}
\n Denote by $\tilde\psi$ the map $\tilde\psi=\psi\compo f_+:
\f_+/\i_+\to\udjp$.
\vskip .3 truecm
\rem{\pshom}{}
\n $\tilde\psi$ is obviously homogeneous, that is $\tilde\psi=\oplus_{\alpha\in Q_+}\tilde\psi_{\alpha}$ with $\tilde\psi_{\alpha}=\tilde\psi\big|_{(\f_+/\i_+)_{\alpha}}$ and consequently
$\tilde\psi_1=\oplus_{\alpha\in Q_+}(\tilde\psi_{\alpha})_1$ where $(\tilde\psi_{\alpha})_1$ is the specialization at 1 of $\tilde\psi_{\alpha}$.

\n Since $(\f_+/\i_+)_{\alpha}$ is finitely generated over $\a$ and $\u_{\a,\alpha}^{DJ,+}$ is free over $\a$ we have that for each $\alpha\in Q_+$ $\tilde\psi_{\alpha}$ is injective if 
$(\tilde\psi_{\alpha})_1$ is injective (see proposition \specgen).

\n Then $\tilde\psi$ is injective if 
$\tilde\psi_1$ is injective.

\vskip .5 truecm

\noindent {\bf{\spquno.\ {\bf SPECIALIZATION at $q=1$.}}}
\vskip .5truecm
\n We are reduced to study the specialization at 1 of $\tilde\psi$. To this aim it is important that first of all we understand the structure of the specialization at 1 of $\f_+/\i_+$ and of $\udjp$. Since, as recalled in remark  \imredd\ below, the specialization at 1 of $\udjp$ is well  known,  we concentrate in the description of the specialization of $\f_+/\i_+$. 

\n Of course a first presentation by generators and relations of the specialization at 1 of $\f_+/\i_+$ is immediate to find by specializing at 1 the defining relations of $\f_+/\i_+$ (see proposition \imredt).
The present section is 
 devoted to simplify these specialized relations.

\vskip .3 truecm
\rem{\imredd}{}
\n Thanks to remark \qremd, iv), the specialization at 1 of $\udjp$ is the enveloping algebra of the Lie algebra generated by $\{e_i|i\in I\}$ with relations $({\rm{ad}}e_i)^{1-a_{ij}}(e_j)=0$ 
when $i\neq j$ (Serre relations),
which is well known to be the positive part of the Kac-Moody algebra $\hat\gothg=\hat\gothg(X_{\tilde n}^{(k)})$ and also of the loop algebra $(\gothg\otimes_{\C}\C[t^{\pm 1}])^{\chi}\supseteq\gothg^{\chi}=\gothg_0$ (see \kac).
\vskip .3 truecm
\prop{\imredt}
\n By the very definition of $\f_+$ and $\i_+$ the specialization at 1 of $\f_+/\i_+$ is the (associative) algebra 
generated by $\{x_{i,r}^+|i\in I_0, r\geq 0\}$ with the following relations $({\tilde{dr}}_+)$:
$$x_{i,r}^+=0\ (\tilde d_i\not|r),\leqno{(zx)}$$
$$[x_{i,r+\tilde d_{ij}}^{+},x_{j,s}^{+}]
+[x_{j,s+\tilde d_{ij}}^{+},x_{i,r}^{+}]=0\ ((i,r),(j,s)\in \iz,\ 
a_{ij}\!\!<\!\!0),\leqno{(xd^{})}$$
$$
\sum_{\sigma\in\sy_2}\sigma.[x_{i,r_1+\tilde d_i}^{+},x_{i,r_2}^{+}]=0\ \ ((r_1,r_2)\in\N^2,\ (X_{\tilde n}^{(k)},d_i)\neq(A_{2n}^{(2)},1)),\leqno{(x1^{})}$$
$$
\sum_{\sigma\in\sy_2}\sigma.[x_{i,r_1+ 2}^{+},x_{i,r_2}^{+}]=0\ \ \ \ ((r_1,r_2)\in\N^2,\ (X_{\tilde n}^{(k)},d_i)=(A_{2n}^{(2)},1)),\leqno{(x2^{})}$$
$$\sum_{\sigma\in\sy_3}\sigma.[[x_{i,r_1+ 1}^{+},x_{i,r_2}^{+}],x_{i,r_3}^{+}]=0\leqno{(x3)}$$
$\ \ \ \ ((r_1,r_2,r_3),\in\N^3,\ (X_{\tilde n}^{(k)},d_i)=(A_{2n}^{(2)},1)),$
$$\sum_{\sigma\in\sy_2}\sigma.[[x_{j,s}^{+},x_{i,r_1+ 1}^{+}],x_{i,r_2}^{+}]=0,\leqno(t2^{})$$
$$\ \ \ \ (i,j\in I_0,\ a_{ij}=-2,\ k=2,\ X_{\tilde n}^{(k)}\neq A_{2n}^{(2)},\ (r_1,r_2)\in\N^{2},\ s\in\Z),$$
$$\sum_{\sigma\in\sy_2}\sigma.\big([[x_{j,s}^{+},x_{i,r_1+ 1}^{+}],x_{i,r_2}^{+}]+[[x_{j,s}^{+},x_{i,r_2}^
{+}],x_{i,r_1+1}^{+}])=0\leqno(s2^{})$$
$$\ \ \ \ (i,j\in I_0,\ a_{ij}=-2,\ X_{\tilde n}^{(k)}= A_{2n}^{(2)},\ (r_1,r_2)\in\N^{2},\ s\in\Z),$$
$$\sum_{\sigma\in\sy_2}\sigma.(2[[x_{j,s}^{+},x_{i,r_{1}+ 2}^{+}],x_{i,r_{2}}^{+}]+[[x_{j,s}^{+},x_{i,r_{1}+ 1}^{+}],x_{i,r_{2}+ 1}^{+}])=0\leqno{(t3^{})}$$
$$\ \ \ \ (i,j\in I_0,\ a_{ij}=-3,\ k=3,\ (r_1,r_2)\in\N^{2},\ s\in\Z),$$
$$\sum_{\sigma\in\sy_3}\sigma.[[[x_{j,s}^{+},x_{i,r_1+ 1}^{+}],x_{i,r_2}^{+}],x_{i,r_3}^{+}]=0,\leqno(u3^{})$$
$$\ \ \ \ (i,j\in I_0,\ a_{ij}=-3,\ k=3,\ (r_1,r_2,r_3)\in\N^{3},\ s\in\Z),$$
$$\sum_{\sigma\in\sy_{1-a_{ij}}}\sigma.
[...[[x_{j,s}^+,x_{i,r_1}^{+}],x_{i,r_2}^+],...
,x_{i,r_{1-a_{ij}}}^{+}]
=0\leqno{(s)}$$
$$\ \ \ \ (i\neq j\in\! I_0,
\ r\in
\N^{1-a_{ij}},\ s\in\Z).$$
\dim
All the relations $({\tilde{dr}}_+)$ are the immediate specialization at 1 of the relations $(ZX_+^+, DR_+,U3_+,S_+)$, recalling notation \nbr\ and remark \dremd,ii), and noticing that relations $(S2_+)$ specialize to
$$0=\sum_{\sigma\in\sy_2}\Big(2[[x_{j,s}^+,x_{i,r_1+1}^+],x_{i,r_2}^+]+[[x_{i,r_1+1}^+,x_{i,r_2}^+],x_{j,s}^+]\Big)=$$
$$=\sum_{\sigma\in\sy_2}\Big(2[[x_{j,s}^+,x_{i,r_1+1}^+],x_{i,r_2}^+]+[x_{i,r_1+1}^+,[x_{i,r_2}^+,x_{j,s}^+]]-[x_{i,r_2}^+,[x_{i,r_1+1}^+,x_{j,s}^+]]\Big)=$$
$$=\sum_{\sigma\in\sy_2}\Big([[x_{j,s}^+,x_{i,r_1+1}^+],x_{i,r_2}^+]+[x_{i,r_1+1}^+,[x_{i,r_2}^+,x_{j,s}^+]]\Big),$$
which is $(s2)$.

\vskip .3 truecm
\rem{\imredq}{}
\n Remark that in the relations $({\tilde {dr}}_+)$ (see proposition \imredt) all the products are expressed in terms of brackets; hence the 
associative algebra generated by $\{x_{i,r}^+|i\in I_0, r\geq 0\}$ with the relations $({\tilde{dr}}_+)$
is the enveloping algebra of the Lie algebra 
generated by $\{x_{i,r}^+|i\in I_0, r\geq 0\}$ with the relations $({\tilde{dr}}_+)$.

\n This Lie algebra plays a central role in the following.
\vskip .3 truecm
\ddefi{\aggde}
\n $L_+$ is the Lie algebra generated by $\{x_{i,r}^+|i\in I_0, r\geq 0\}$ with the relations $({\tilde{dr}}_+)$.
\vskip .3 truecm
\rem{\aggre}{}
\n The specialization at 1 of $\f_+/\i_+$ is the enveloping algebra $\u(L_+)$ of the Lie 
algebra $L_+$ (see proposition \imredt, remark \imredq\ and definition \aggde). 

\n In particular $\tilde\psi_1$ is a homomorphism of associative algebras from $\u(L_+)$ to $\u(\hat\gothg^+)$, see remark \imredd.
The next step is proving that $\tilde\psi_1(L_+)\subseteq\hat\gothg^+$, which implies that $\tilde\psi_1\big|_{L_+}$  is a Lie algebra homomorphism from $L_+$ to $\hat\gothg^+$ and 
$\tilde\psi_1=\u(\tilde\psi_1\big|_{L_+})$.

\vskip .3 truecm
\rem{\imredc}{}
\n $\tilde\psi_1(L_+)\subseteq(\gothg_+\otimes_{\C}\C[t])^{\chi}\subseteq\hat\gothg^+$; in particular, thanks to remark \aggre\ we obtain that $\tilde\psi_1$ is injective if and only if $\tilde\psi_1\big|_{L_+}$ is injective.
\dim
Since $\tilde\psi\compo\bar t_+^{{\prime}}=T_{\lambda}^{-1}\compo\tilde\psi$ (see remarks \repsi, iii) and \firem,iv)),
the claim follows from the fact that $\tilde\psi(x_{i,0}^+)=e_i\in\gothg_{0,+}\subseteq\hat\gothg$, from remarks \qremd,viii) and \drromrem,ix), and from the fact that $\hat\gothg_{r\delta+\alpha_i}\subseteq(\gothg_+\otimes\C[t])^{\chi}$ if $i\in I_0$, $r\geq 0$.

\vskip .3 truecm
\prop{\imreds}
\n $\tilde\psi_1\big|_{L_+}:L_+\to(\gothg_+\otimes_{\C}\C[t])^{\chi}$ is surjective.
\dim
$(\gothg_+\otimes_{\C}\C[t])^{\chi}=\oplus_{r\in\N}(\gothg_+^{[r]}\otimes_{\C}\C t^r)\subseteq\oplus_{r\in\N}(\gothg^{[r]}\otimes_{\C}\C t^r)=(\gothg\otimes_{\C}\C[t])^{\chi}$
where $\gothg^{[r]}=\gothg_-^{[r]}\oplus\gothh^{[r]}\oplus\gothg_+^{[r]}$ is well known to be a simple finite dimensional $\gothg_0=\gothg^{[0]}$-module, hence a lowest weight cyclic $\gothg_+^{[0]}=\gothg_{0,+}$-module (see \kac).

\n Then $\gothg_+^{[r]}(=\oplus_{\alpha\in Q_{0,+}\setminus\{0\}}(\gothg^{[r]})_{\alpha})$
is generated as a $\gothg_{0,+}$-module by $$\oplus_{i\in I_0}(\gothg^{[r]})_{\alpha_i}=\oplus_{i\in I_0}(\gothg_+^{[r]})_{\alpha_i}\big(=\oplus_{{i\in I_0:\atop \tilde d_i|r}}(\gothg_+^{[r]})_{\alpha_i}\ {\rm{since\ }}(\gothg_+^{[r]})_{\alpha_i}=(0)\ {\rm{if}}\ \tilde d_i\not| r\big),$$
that is $(\gothg_+\otimes_{\C}\C[t])^{\chi}$ 
is generated as a $\gothg_{0,+}$-module by
$\oplus_{{i\in I_0,r\in\N:\atop \tilde d_i|r}}(\gothg_+^{[r]})_{\alpha_i}\otimes\C t^r$ or equivalently by
$\{\tilde\psi_1(x_{i,r}^+)|i\in I_0,r\in\N$ such that $\tilde d_i|r\}$ since $\forall i\in I_0, r\in\N$ $\tilde\psi_1(x_{i,\tilde d_ir}^+)=\pm\tilde T_{\lambda}^{-r}(e_i)\neq 0$ 
(see remarks \qremd,viii), \drromrem,ix), \repsi, iii) and \firem,iv)),
and $(\gothg_+^{[\tilde d_ir]})_{\alpha_i}$ is one dimensional.

\n This forces $\{\tilde\psi_1(x_{i,r}^+)|i\in I_0,r\in\N\}$, which obviously contains $\{e_i=\tilde\psi_1
(x_{i,0}^+)|i\in I_0\}$, to generate $(\gothg_+\otimes_{\C}\C[t])^{\chi}$ also as a Lie algebra; the assertion follows.
\vskip .3 truecm
\cor{\epip}
\n i) $E_{r\delta+\alpha}\in\psi(\uq^{Dr,+,+})$ if $r\geq 0$ and $\alpha\in Q_{0,+}\setminus\{0\}$;

\n ii) $\uq^{DJ,+}\cap\psi(\uq^{Dr,+})=\psi(\uq^{Dr,+,+})$.
\dim
i) follows from ii) by corollary \pip\ (indeed i) and ii) are equivalent claims because $\psi(\uq^{Dr,+,+})\subseteq\uq^{DJ,+}\cap\psi(\uq^{Dr,+})$). So it is enough to compare the dimensions of the homogeneous components of $\uq^{DJ,+}\cap\psi(\uq^{Dr,+})$ and $\psi(\uq^{Dr,+,+})$: for all $\alpha\in Q$
$${\rm{dim}}_{\C(q)}\u_{q,\alpha}^{DJ,+}\cap\psi(\uq^{Dr,+})\geq{\rm{dim}}_{\C(q)}\psi(\u_{q,\alpha}^{Dr,+,+})=rk_{\a}\tilde\psi(\f_+/\i_+)_{\alpha}=$$
$$={\rm{dim}}_{\C}\tilde\psi_1(\u(L_+)_{\alpha})={\rm{dim}}_{\C}\u((\gothg_+\otimes_{\C}\C[t])^{\chi})_{\alpha}={\rm{dim}}_{\C(q)}\u_{q,\alpha}^{DJ,+}\cap\psi(\uq^{Dr,+})$$
where the last two equalities follow respectively from proposition \imreds\ and from the comparison of  the $\C(q)$-basis of $\u_q^{DJ,+}\cap\psi(\uq^{Dr,+})$ described in corollary \pip\ with the PBW-basis of 
$\u((\gothg_+\otimes_{\C}\C[t])^{\chi})$.

\vskip .3 truecm
\n Before proving, in section \affkm, that  the Lie-algebra homomorphism $\tilde\psi_1\big|_{L_+}$ is actually injective, in the remaining part of this section we simplify the relations defining $L_+$ (see the following computations, summarized in corollary \lpgen).
\vskip .3 truecm
\rem{\imredo}{}
\n Relations $(xd)$ are equivalent to saying that if $a_{ij}<0$, $\tilde d_i|r$ and $\tilde d_j|s$ ($i\neq j$ fixed) then $[x_{i,r}^+,x_{j,s}^+]$ depends only on $r+s$.
Together with $(s)$ in case $a_{ij}=0$ they imply 
$$[x_{i,r}^+,x_{j,s}^+]\ {\rm{depends\ only\ on\ }}r+s\ \ (i\neq j\in I_0\ {\rm {fixed}},\ \tilde d_i|r,\ \tilde d_j|s).\leqno{(x_d)}$$
\vskip .3 truecm
\lem{\imleu}
\n Relations $(x1)$ and $(x2)$ are equivalent to
$$
[x_{i,r}^{+},x_{i,s}^{+}]\!=\!\begin{cases}\!0&\!\!\!\!{\rm{if}}\ (X_{\tilde n}^{(k)},d_i)\neq(A_{2n}^{(2)},1)\ {\rm{or}}\ 2|r+s\cr \!(-1)^h[x_{i,s+h+1}^{+},x_{i,s+h}^{+}]&\!\!\!\!
{\rm{if
}}\ r=s+2h+1;
\end{cases}\leqno{(x_{1,2})}$$
in particular $(-1)^s[x_{i,r}^{+},x_{i,s}^{+}]
$ depends only on $r+s$.
\dim
That $(x_{1,2})$ implies $(x1)$ and $(x2)$ is obvious. Viceversa:

\n i) case $(X_{\tilde n}^{(k)},d_i)\neq(A_{2n}^{(2)},1)$: of course we can suppose $r\geq s$ and proceed by induction on $r-s$, the cases $r=s$ and $r=s+\tilde d_i$ being obvious; if $r>s+\tilde d_i$
$$[x_{i,r}^+,x_{i,s}^+]=-[x_{i,s+\tilde d_i}^+,x_{i,r-\tilde d_i}^+]=[x_{i,r-\tilde d_i}^+,x_{i,s+\tilde d_i}^+]=0\ \ (r-\tilde d_i\geq s+\tilde d_i).$$ 
\n ii) case $(X_{\tilde n}^{(k)},d_i)=(A_{2n}^{(2)},1)$: again we can suppose $r\geq s$ and proceed by induction on $r-s$, the cases $r-s=0,1,2$ being obvious:
$$r-s=3\Rightarrow[x_{i,r}^+,x_{i,s}^+]=[x_{i,s+3}^+,x_{i,s}^+]=-[x_{i,s+2}^+,x_{i,s+1}^+];$$
$r-s>3\Rightarrow[x_{i,r}^+,x_{i,s}^+]=-[x_{i,s+2}^+,x_{i,r-2}^+]=[x_{i,r-2}^+,x_{i,s+2}^+]$, from which the claim follows by the inductive hypothesis, since $r-2\geq s+2$.
\vskip .3 truecm
\cor{\imcou}
\n If $(X_{\tilde n}^{(k)},d_i)\neq(A_{2n}^{(2)},1)$ or $2|r+s$ we have 
$[[a,x_{i,r}^+],x_{i,s}^+]=[[a,x_{i,s}^+],x_{i,r}^+]$ for all $a\in L_+$.
\dim
Indeed $[[a,x_{i,r}^+],x_{i,s}^+]-[[a,x_{i,s}^+],x_{i,r}^+]=[a,[x_{i,r}^+,x_{i,s}^+]]=0$ thanks to lemma \imleu.

\vskip .3 truecm
\lem{\imled}
\n If $(X_{\tilde n}^{(k)},d_i)=(A_{2n}^{(2)},1)$ relations $(x_{1,2})$, $(x3)$ imply
$$[[x_{i,r_1}^{+},x_{i,r_2}^{+}],x_{i,r_3}^{+}]=0\ \ ((r_1,r_2,r_3)\in\N^3
).\leqno{(x_3)}$$
\dim
Thanks to lemma \imleu\ it is enough to prove that 
$$[[x_{i,r+1}^+,x_{i,r}^+],x_{i,s}^+]=0 \ \forall r,s\in\N.$$
Recall that by $(x3)$
$$[[x_{i,r+1}^+,x_{i,r}^+],x_{i,s}^+]+[[x_{i,r+1}^+,x_{i,s}^+],x_{i,r}^+]+[[x_{i,s+1}^+,x_{i,r}^+],x_{i,r}^+]=0;$$
if $r+s+1$ is even then, by $(x_{1,2})$,
$$[x_{i,r+1}^+,x_{i,s}^+]=0=[x_{i,s+1}^+,x_{i,r}^+]$$ so that 
$$[[x_{i,r+1}^+,x_{i,r}^+],x_{i,s}^+]=0;$$ if $r+s$ is even then by corollary \imcou
$$[[x_{i,r+1}^+,x_{i,r}^+],x_{i,s}^+]=[[x_{i,r+1}^+,x_{i,s}^+],x_{i,r}^+];$$
moreover by  $(x_{1,2})$ $[x_{i,s+1}^+,x_{i,r}^+]=\pm[x_{i,r+1}^+,x_{i,s}^+]$,
so that 
$$0=[[x_{i,r+1}^+,x_{i,r}^+],x_{i,s}^+]+[[x_{i,r+1}^+,x_{i,s}^+],x_{i,r}^+]+[[x_{i,s+1}^+,x_{i,r}^+],x_{i,r}^+]=$$
$$=(2\pm 1))[[x_{i,r+1}^+,x_{i,s}^+],x_{i,r}^+],$$
which is $[[x_{i,r+1}^+,x_{i,r}^+],x_{i,s}^+]=0$.
\vskip .3 truecm
\prop{\imredn}
\n Relations $(x1)$, $(x2)$, $(x3)$ are equivalent to relations $(x_{1,2})$, $(x_3)$ (it is obvious that $(x_3)$ implies $(x3)$).

\vskip .3 truecm
\lem{\emo}
\n Let $i,j\in I_0$, $r_1,r_2,s\in\Z$ be such that $a_{ij}<0$, $\tilde d_i|r_1,r_2$ and $\tilde d_j|s$. Then:

\n i) if $\tilde d_i\geq\tilde d_j$ and $(X_{\tilde n}^{(k)},d_i)\neq(A_{2n}^{(2)},1)$
$$[[x_{j,s}^+,x_{i,r_1}^+],x_{i,r_2}^+]=[[x_{j,s+r_1+r_2}^+,x_{i,0}^+],x_{i,0}^+];$$
ii)  if $1=\tilde d_i<\tilde d_j=k$ or $(X_{\tilde n}^{(k)},d_i)=(A_{2n}^{(2)},1)$, and $k|r_2-\varepsilon_2$ ($0\leq\varepsilon_2<k$) then
$[[x_{j,s}^+,x_{i,r_1}^+],x_{i,r_2}^+]$ depends only on $(s+r_1+r_2,\varepsilon_2)$.

\dim
i) is an immediate consequence of relations $(x_d)$, $(x_{1,2})$ and of corollary \imcou:
$$[[x_{j,s}^+,x_{i,r_1}^+],x_{i,r_2}^+]=[[x_{j,s+r_1}^+,x_{i,0}^+],x_{i,r_2}^+]=$$
$$=[[x_{j,s+r_1}^+,x_{i,r_2}^+],x_{i,0}^+]=[[x_{j,s+r_1+r_2}^+,x_{i,0}^+],x_{i,0}^+];$$
ii) is similar: if $\tilde d_i<\tilde d_j$ or $2|r_1+r_2$
$$[[x_{j,s}^+,x_{i,r_1}^+],x_{i,r_2}^+]=[[x_{j,s}^+,x_{i,r_2}^+],x_{i,r_1}^+]=[[x_{j,s+r_2-\varepsilon_2}^+,x_{i,\varepsilon_2}^+],x_{i,r_1}^+]=$$
$$=[[x_{j,s+r_2-\varepsilon_2}^+,x_{i,r_1}^+],x_{i,\varepsilon_2}^+];$$
if $(X_{\tilde n}^{(k)},d_i)=(A_{2n}^{(2)},1)$, and $s+r_1>0$ or $2|r_2$
$$[[x_{j,s}^+,x_{i,r_1}^+],x_{i,r_2}^+]=[[x_{j,s+r_1-\varepsilon_2}^+,x_{i,\varepsilon_2}^+],x_{i,r_2}^+]=[[x_{j,s+r_1-\varepsilon_2}^+,x_{i,r_2}^+],x_{i,\varepsilon_2}^+];$$
in both cases $[[x_{j,s}^+,x_{i,r_1}^+],x_{i,r_2}^+]$ depends only on $(s+r_1+r_2-\varepsilon_2,\varepsilon_2)$, that is on $(s+r_1+r_2,\varepsilon_2)$;

\n finally if $(X_{\tilde n}^{(k)},d_i)=(A_{2n}^{(2)},1)$ and 
$s=r_1=0$, $r_2=2r+1$ we can suppose $r>0$ and we have
$$[[x_{j,0}^+,x_{i,0}^+],x_{i,2r+1}^+]=[[x_{j,0}^+,x_{i,2r+1}^+],x_{i,0}^+]+[x_{j,0}^+,[x_{i,0}^+,x_{i,2r+1}^+]]=$$
$$=[[x_{j,0}^+,x_{i,2r+1}^+],x_{i,0}^+]+[x_{j,0}^+,[x_{i,2}^+,x_{i,2r-1}^+]]=$$
$$=
[[x_{j,0}^+,x_{i,2r+1}^+],x_{i,0}^+]+[[x_{j,0}^+,x_{i,2}^+],x_{i,2r-1}^+]-[[x_{j,0}^+,x_{i,2r-1}^+],x_{i,2}^+]=$$
$$=[[x_{j,0}^+,x_{i,2}^+],x_{i,2r-1}^+]$$
and the claim follows from the previous cases.

\vskip .3 truecm
\prop{\emu}
\n Relations $(x_d)$, $(x_{1,2})$ and $(t2)$ are equivalent to relations $(x_d)$, $(x_{1,2})$,  $(t_2^{\prime})$ and $(t_2^{\prime\prime})$, where $(t_2^{\prime})$ and $(t_2^{\prime\prime})$ are the following relations:
$$[[x_{j,s}^+,x_{i,1}^+],x_{i,0}^+]=0
,\leqno{(t_2^{\prime}
)}$$
$$[[x_{j,s}^+,x_{i,1}^+],x_{i,1}^+]=-[[x_{j,s+2}^+,x_{i,0}^+],x_{i,0}^+]\leqno{(t_2^{\prime\prime}
)}$$
$(k=2,\ a_{ij}=-2, X_{\tilde n}^{(k)}\neq A_{2n}^{(2)})$.
\dim
Indeed corollary \imcou\ and remark \emo\ imply that 
$$\sum_{\sigma\in\sy_2}\sigma.[[x_{j,s}^+,x_{i,r_1+1}^+],x_{i,r_2}^+]=[[x_{j,s}^+,x_{i,r_1+1}^+],x_{i,r_2}^+]+[[x_{j,s}^+,x_{i,r_1}^+],x_{i,r_2+1}^+]=$$
$$=\begin{cases}[[x_{j,s+r_1+r_2}^+,x_{i,1}^+],x_{i,0}^+]+[[x_{j,s+r_1+r_2}^+,x_{i,0}^+],x_{i,1}^+]&{\rm{if}}\ 2|r_1+r_2\cr
[[x_{j,s+r_1+r_2-1}^+,x_{i,1}^+],x_{i,1}^+]+[[x_{j,s+r_1+r_2+1}^+,x_{i,0}^+],x_{i,0}^+]&{\rm{otherwise;}}
\end{cases}$$
but by corollary \imcou\ we have
$$[[x_{j,s+r_1+r_2}^+,x_{i,1}^+],x_{i,0}^+]+[[x_{j,s+r_1+r_2}^+,x_{i,0}^+],x_{i,1}^+]=2[[x_{j,s+r_1+r_2}^+,x_{i,1}^+],x_{i,0}^+].$$

\vskip .3 truecm
\prop{\emd}
\n Relations $(x_d)$, $(x_{1,2})$ and $(s2)$ are equivalent to relations $(x_d)$, $(x_{1,2})$ and $(s_2)$, where $(s_2)$ are the following relations:
$$[[x_{j,s}^+,x_{i,0}^+],x_{i,1}^+]+[[x_{j,s+1}^+,x_{i,0}^+],x_{i,0}^+]=0\ \ 
\ \ (a_{ij}=-2, 
X_{\tilde n}^{(k)}= A_{2n}^{(2)})
;\leqno{(s_2
)}$$
\dim
Indeed lemma \emo,ii) implies that 
$$\sum_{\sigma\in\sy_2}\sigma.([[x_{j,s}^+,x_{i,r_1+1}^+],x_{i,r_2}^+]+[[x_{j,s}^+,x_{i,r_2}^+],x_{i,r_1+1}^+])=$$
$$=2([[x_{j,s+r_1+r_2}^+,x_{i,0}^+],x_{i,1}^+]+[[x_{j,s+r_1+r_2+1}^+,x_{i,0}^+],x_{i,0}^+]).$$

\vskip .3 truecm
\lem{\nqud}
\n Relations $(x_d)$, $(x_{1,2})$ and $(t3)$ are equivalent to relations $(x_d)$, $(x_{1,2})$ and $(\tilde t3)$ where relations $(\tilde t3)$ are the following:
$$[[x_{j,s}^+,x_{i,r_1+2}^+],x_{i,r_2}^+]+[[x_{j,s}^+,x_{i,r_1+1}^+],x_{i,r_2+1}^+]+[[x_{j,s}^+,x_{i,r_1}^+],x_{i,r_2+2}^+]=0\leqno{(\tilde t3)}$$
($k=3$, $a_{ij}=-3$, $s,r_1,r_2,r_3\in\N$).
\dim
Indeed by corollary \imcou\
$$\sum_{\sigma\in\sy_2}\sigma.(2[[x_{j,s}^+,x_{i,r_1+2}^+],x_{i,r_2}^+]+[[x_{j,s}^+,x_{i,r_1+1}^+],x_{i,r_2+1}^+])=$$
$$=2([[x_{j,s}^+,x_{i,r_1+2}^+],x_{i,r_2}^+]+[[x_{j,s}^+,x_{i,r_1}^+],x_{i,r_2+2}^+]+[[x_{j,s}^+,x_{i,r_1+1}^+],x_{i,r_2+1}^+]).$$

\vskip .3 truecm
\nota{\ntnz}{}
\n Let us define the following relations:
$$[[x_{j,s}^+,x_{i,1}^+],x_{i,1}^+]=-2[[x_{j,s}^+,x_{i,2}^+],x_{i,0}^+]\ \ (k=3,\ a_{ij}=-3)\leqno{(t_3^{\prime})}$$
$$2[[x_{j,s}^+,x_{i,2}^+],x_{i,1}^+]=-[[x_{j,s+3}^+,x_{i,0}^+],x_{i,0}^+]\ \ (k=3,\ a_{ij}=-3)\leqno{(t_3^{\prime\prime})}$$
$$[[x_{j,s}^+,x_{i,2}^+],x_{i,2}^+]=-2[[x_{j,s+3}^+,x_{i,1}^+],x_{i,0}^+]\ \ (k=3,\ a_{ij}=-3)\leqno{(t_3^{\prime\prime\prime})}$$

\vskip .3 truecm
\rem{\relsem}{}
\n Relations $(x_d)$, $(x_{1,2})$, $(\tilde t3)$ imply relations $(t_3^{\prime})$-$(t_3^{\prime\prime\prime})$.
\dim
\n Using relations $(x_d)$ and $(x_{1,2})$ we have of course that 
$(t_3^{\prime})$, $(t_3^{\prime\prime})$ and $(t_3^{\prime\prime\prime})$ are $(\tilde t3)$ with $r_1+r_2=0,1,2$ respectively.

\vskip .3 truecm
\prop{\nrku}
\n Relations $(x_d)$, $(x_{1,2})$, $(t_3^{\prime})$, $(t_3^{\prime\prime})$,  $(t_3^{\prime\prime\prime})$ are equivalent to relations $(x_d)$, $(x_{1,2})$, $(t3)$.
\dim
\n We prove by induction on $r_1+r_2$ that relations $(x_d)$, $(x_{1,2})$, $(t_3^{\prime})$, $(t_3^{\prime\prime})$,  $(t_3^{\prime\prime\prime})$ imply relations $(\tilde t3)$, the cases $0\leq r_1+r_2<3$ being obvious (see the proof of remark \relsem).
If $r_1+r_2\geq 3$ use induction on $r_2$: if $r_2=0$ then $r_1\geq 3$ and thanks to $(x_d)$ we have
$$[[x_{j,s}^+,x_{i,r_1+2}^+],x_{i,r_2}^+]+[[x_{j,s}^+,x_{i,r_1+1}^+],x_{i,r_2+1}^+]+[[x_{j,s}^+,x_{i,r_1}^+],x_{i,r_2+2}^+]=$$
$$=[[x_{j,s+3}^+,x_{i,r_1-1}^+],x_{i,0}^+]+[[x_{j,s+3}^+,x_{i,r_1-2}^+],x_{i,1}^+]+[[x_{j,s+3}^+,x_{i,r_1-3}^+],x_{i,2}^+],$$
which is zero by the inductive hypothesis ($r_1-3+0<r_1+r_2)$;
if $r_2>0$ then,  thanks to lemma \emo, iii), 
$[[x_{j,s}^+,x_{i,r_1}^+],x_{i,r_2+2}^+]=
[[x_{j,s}^+,x_{i,r_1+3}^+],x_{i,r_2-1}^+]$, so that
$$[[x_{j,s}^+,x_{i,r_1+2}^+],x_{i,r_2}^+]+[[x_{j,s}^+,x_{i,r_1+1}^+],x_{i,r_2+1}^+]+[[x_{j,s}^+,x_{i,r_1}^+],x_{i,r_2+2}^+]=$$
$$=[[x_{j,s}^+,x_{i,r_1+3}^+],x_{i,r_2-1}^+]+[[x_{j,s}^+,x_{i,r_1+2}^+],x_{i,r_2}^+]+[[x_{j,s}^+,x_{i,r_1+1}^+],x_{i,r_2+1}^+]$$
which is zero because $r_2-1<r_2$.
\vskip .3 truecm
\rem{\nlmu}{}
\n If $k=3$, $a_{ij}=-3$ relations $(x_d)$, $(x_{1,2})$, $(t_3^{\prime})$, $(t_3^{\prime\prime})$, $(t_3^{\prime\prime\prime})$ imply that if $s>0$ then
$$[[[x_{j,s}^+,x_{i,1}^+],x_{i,0}^+],x_{i,0}^+]=0
$$
and $$[[[x_{j,s}^+,x_{i,2}^+],x_{i,0}^+],x_{i,0}^+]=0.
$$
\dim
$s>0\Rightarrow s\geq 3$, then relations $(t_3^{\prime\prime\prime})$, $(t_3^{\prime})$, $(t_3^{\prime\prime})$ imply that 
$$2[[[x_{j,s}^+,x_{i,1}^+],x_{i,0}^+],x_{i,0}^+]+[[[x_{j,s-3}^+,x_{i,2}^+],x_{i,2}^+],x_{i,0}^+]=0$$
$$2[[[x_{j,s-3}^+,x_{i,2}^+],x_{i,0}^+],x_{i,2}^+]+[[[x_{j,s-3}^+,x_{i,1}^+],x_{i,1}^+],x_{i,2}^+]=0$$
$$2[[[x_{j,s-3}^+,x_{i,2}^+],x_{i,1}^+],x_{i,1}^+]+[[[x_{j,s}^+,x_{i,0}^+],x_{i,0}^+],x_{i,1}^+]=0,$$
from which, thanks to corollary \imcou, 
$$9[[[x_{j,s}^+,x_{i,1}^+],x_{i,0}^+],x_{i,0}^+]=0.$$
Analogously 
$$2[[[x_{j,s}^+,x_{i,2}^+],x_{i,0}^+],x_{i,0}^+]+[[[x_{j,s}^+,x_{i,1}^+],x_{i,1}^+],x_{i,0}^+]=0$$
$$2[[[x_{j,s}^+,x_{i,1}^+],x_{i,0}^+],x_{i,1}^+]+[[[x_{j,s-3}^+,x_{i,2}^+],x_{i,2}^+],x_{i,1}^+]=0$$
$$2[[[x_{j,s-3}^+,x_{i,2}^+],x_{i,1}^+],x_{i,2}^+]+[[[x_{j,s}^+,x_{i,0}^+],x_{i,0}^+],x_{i,2}^+]=0,$$
from which $[[[x_{j,s}^+,x_{i,2}^+],x_{i,0}^+],x_{i,0}^+]=0.$
\vskip .3 truecm

\nota{\ntzu}{}
\n Let us define the following relations:
$$[[[x_{j,0}^+,x_{i,1}^+],x_{i,0}^+],x_{i,0}^+]=0\ \ (k=3,\ a_{ij}=-3),\leqno{(u_3^{\prime})}$$
$$[[[x_{j,0}^+,x_{i,2}^+],x_{i,0}^+],x_{i,0}^+]=0\ \ (k=3,\ a_{ij}=-3).\leqno{(u_3^{\prime\prime})}
$$
\vskip .3 truecm
\rem{\nrkt}{}
\n Relations $(x_{1,2})$, $(t_3^{\prime})$ and $(u3)$ imply relations $(u_3^{\prime})$ and $(u_3^{\prime\prime})$.
\dim
\n $(u_3^{\prime})$ is $(u3)$ with $s=r_1=r_2=r_3=0$.

\n $[(t_3^{\prime}),x_{i,0}^+]$ with $s=0$ and  $(u3)$ with $(s,r_1,r_2,r_3)=(0,1,0,0)$ imply $(u_3^{\prime\prime})$, using corollary \imcou.

\vskip .3 truecm
\prop{\nrkq}
\n Relations $(x_d)$, $(x_{1,2})$, $(t_3^{\prime})$, $(t_3^{\prime\prime})$,  $(t_3^{\prime\prime\prime})$, $(u_3^{\prime})$, $(u_3^{\prime\prime})$ imply relations $(u3)$ (hence are equivalent to relations $(x_d)$, $(x_{1,2})$, $(t3)$, $(u3)$).
\dim
The hypotheses imply that $[[[x_{j,s}^+,x_{i,r_1}^+],x_{i,r_2}^+],x_{i,r_3}^+]$ is a rational multiple of $[[[x_{j,s}^+,x_{i,r_1+r_2+r_3}^+],x_{i,0}^+],x_{i,0}^+]$ (by $(x_d)$, $(x_{1,2})$, $(t_3^{\prime})$, $(t_3^{\prime\prime})$,  $(t_3^{\prime\prime\prime})$), which is zero if $3\not|r_1+r_2+r_3$ (by $(x_d)$, $(u_3^{\prime})$, $(u_3^{\prime\prime})$ and lemma \nlmu). In particular $(u3)$ holds if 
$3\not|r_1+r_2+r_3+1$. Otherwise we can suppose $0\leq r_1,r_2, r_3<3$, $r_1=r_2$, $r_3+1\equiv r_1({\rm{mod}}3)$ (thanks to $(x_d)$ and $(x_{1,2})$), or equivalently that $(r_1,r_2,r_3)=(1,1,0),\ (2,0,0),\ (2,2,1)$. In these cases $(u3)$ corresponds respectively to $[(t_3^{\prime}),x_{i,1}^+]$, $[(t_3^{\prime\prime}),x_{i,0}^+]$, and $[(t_3^{\prime\prime\prime}),x_{i,2}^+]$.
\vskip .3 truecm
\nota{\ntzd}{}
\n Let us define the following relations:
$$({\rm {ad}}x_{i,0}^+)^{1-a_{ij}}(x_{j,s}^+)=0\ \ (i\neq j).\leqno{(serre)}$$
\vskip .3 truecm
\rem{\nrkc}{}
\n Relation $(s)$ implies relation $(serre)$.
\dim
\n the claim is obvious since $(serre)$ is $(s)$ with $r_u=0$ for all $u=1,...,1-a_{ij}$.

\vskip .3 truecm
\prop{\imlec}
\n Relations $(x_d)$, $(x_{1,2})$, $(x_{3})$, $(t_2^{\prime})$, $(t_2^{\prime\prime})$, $(s_2)$, $(t_3^{\prime})$, $(t_3^{\prime\prime})$, $(t_3^{\prime\prime\prime})$, $(u_3^{\prime})$, $(u_3^{\prime\prime})$, $(serre)$ imply that $[[...[x_{j,s}^+,x_{i,r_1}^{+}],...],x_{i,r_{1-a_{ij}}}^{+}]=0$ $\forall s\in\Z,r\in\Z^{1-a_{ij}}$; in particular they imply relation $(s)$.
\dim
\n The relations, corollary \imcou\ and lemma \emo\ imply that
$$[[...[x_{j,s}^+,x_{i,r_1}^{+}],...],x_{i,r_{1-a_{ij}}}^{+}]$$ is a rational multiple of $[[...[[x_{j,s}^+,x_{i,r_1+...+r_{1-a_{ij}}}^{+}],x_{i,0}^{+}],...],x_{i,0}^{+}]=$
$$=\begin{cases}({\rm{ad}}x_{i,0}^{+})^{1-a_{ij}}(x_{j,s+r_1+...+r_{1-a_{ij}}}^+)&{\rm{if}}\ \tilde d_j|r_1+...+r_{1-a_{ij}}\cr 0&{\rm{otherwise}}\end{cases}$$ 
(hence zero by $(serre)$) unless in the case $A_{2n}^{(2)}$, $a_{ij}=-2$, $r_3$ odd, when $$[[[x_{j,s}^+,x_{i,r_1}^{+}],x_{i,r_2}^{+}],x_{i,r_{3}}^{+}]=(-1)^{r_2}[[[x_{j,s+r_1+r_2}^+,x_{i,0}^{+}],x_{i,0}^{+}],x_{i,r_{3}}^{+}].$$
But by the above considerations $$[[[x_{j,s}^+,x_{i,0}^{+}],x_{i,0}^{+}],x_{i,r}^{+}]=[[x_{j,s}^+,x_{i,0}^{+}],[x_{i,0}^{+},x_{i,r}^{+}]]+[[[x_{j,s}^+,x_{i,0}^{+}],x_{i,r}^{+}],x_{i,0}^{+}]=$$
$$=[x_{j,s}^+,[x_{i,0}^{+},[x_{i,0}^{+},x_{i,r}^{+}]]]+[[x_{j,s}^+,[x_{i,0}^{+},x_{i,r}^{+}]],x_{i,0}^{+}]=$$
$$=[[[x_{j,s}^+,x_{i,0}^{+}],x_{i,r}^{+}],x_{i,0}^{+}]-[[[x_{j,s}^+,x_{i,r}^{+}],x_{i,0}^{+}],x_{i,0}^{+}]=0.$$
Thus $[[...[x_{j,s}^+,x_{i,r_1}^{+}],...],x_{i,r_{1-a_{ij}}}^{+}]=0$ always.
\vskip .3 truecm
\rem{\nrks}{}
\n It is worth remarking that in the cases $k>1$, $a_{ij}=-k$ relations $(zx)$, $(x_d)$, $(x_{1,2})$, $(x_{3})$, $(t_2^{\prime})$, $(t_2^{\prime\prime})$, $(s_2)$, $(t_3^{\prime})$, $(t_3^{\prime\prime})$, $(t_3^{\prime\prime\prime})$ imply relation $(serre)$ with $s\geq k$ (that is $s>0$ if $X_{\tilde n}^{(k)}\neq A_{2n}^{(2)}$ and $s>1$ if $X_{\tilde n}^{(k)}=A_{2n}^{(2)}$).

\n Compare this observation with
remark \dremd,iii) and iv).
\dim
\n $k=2$, $X_{\tilde n}^{(k)}\neq A_{2n}^{(2)}$: $s>0\Rightarrow s\geq 2$, hence
$$[[[x_{j,s}^+,x_{i,0}^+],x_{i,0}^+],x_{i,0}^+]=-[[[x_{j,s-2}^+,x_{i,1}^+],x_{i,1}^+],x_{i,0}^+]=$$
$$=
-[[[x_{j,s-2}^+,x_{i,1}^+],x_{i,0}^+],x_{i,1}^+]=0;$$
$X_{\tilde n}^{(k)}= A_{2n}^{(2)}$: let $r<s$; then $$[[[x_{j,s}^+,x_{i,0}^+],x_{i,0}^+],x_{i,0}^+]=-[[[x_{j,s-r-1}^+,x_{i,r}^+],x_{i,1}^+],x_{i,0}^+]=$$
$$=-[[[x_{j,s-r-1}^+,x_{i,r}^+],x_{i,0}^+],x_{i,1}^+]-
[[x_{j,s-r-1}^+,x_{i,r}^+],[x_{i,1}^+,x_{i,0}^+]]=$$
$$=-[[[x_{j,s-1}^+,x_{i,0}^+],x_{i,0}^+],x_{i,1}^+]-[[x_{j,s-r-1}^+,[x_{i,1}^+,x_{i,0}^+]],x_{i,r}^+]=$$
$$=-[[[x_{j,s-1}^+,x_{i,0}^+],x_{i,0}^+],x_{i,1}^+]-[[[x_{j,s-r-1}^+,x_{i,1}^+],x_{i,0}^+],x_{i,r}^+]+$$
$$+[[[x_{j,s-r-1}^+,x_{i,0}^+],x_{i,1}^+],x_{i,r}^+]=$$
$$=-[[[x_{j,s-1}^+,x_{i,0}^+],x_{i,0}^+],x_{i,1}^+]-2[[[x_{j,s-r}^+,x_{i,0}^+],x_{i,0}^+],x_{i,r}^+];$$
in particular if $s\geq 2$ we have (choosing $r=0,1$)
$$3[[[x_{j,s}^+,x_{i,0}^+],x_{i,0}^+],x_{i,0}^+]+[[[x_{j,s-1}^+,x_{i,0}^+],x_{i,0}^+],x_{i,1}^+]=0$$
and
$$[[[x_{j,s}^+,x_{i,0}^+],x_{i,0}^+],x_{i,0}^+]+3[[[x_{j,s-1}^+,x_{i,0}^+],x_{i,0}^+],x_{i,1}^+]=0$$
from which the claim follows;

\n $k=3$: $s>0\Rightarrow s\geq 3$, hence
$$[[[x_{j,s}^+,x_{i,0}^+],x_{i,0}^+],x_{i,0}^+],x_{i,0}^+]=-2[[[x_{j,s-3}^+,x_{i,2}^+],x_{i,1}^+],x_{i,0}^+],x_{i,0}^+]=$$
$$=
-2[[[x_{j,s-3}^+,x_{i,2}^+],x_{i,0}^+],x_{i,0}^+],x_{i,1}^+]=0;$$

\vskip .3 truecm
\cor{\lpgen}
\n $L_+$ is the Lie algebra generated by $\{x_{i,r}^+|i\in I_0,\ r\geq 0\}$ with relations $(zx)$, $(x_d)$, $(x_{1,2})$, $(x_3)$, $(t_2^{\prime})$, $(t_2^{\prime\prime})$, $(s_2)$,  $(t_3^{\prime})$,  $(t_3^{\prime\prime})$,  $(t_3^{\prime\prime\prime})$, 
$(u_3^{\prime})$, $(u_3^{\prime\prime})$, $(serre)$.

\vskip .5 truecm
\noindent {\bf{\affkm.\ {\bf AFFINE KAC-MOODY CASE.}}}
\vskip .5truecm

\n This section is devoted to the study of the Lia algebra $L_+$ (see definition \aggde) and of its relation, through $\tilde\psi_1$ (see remarks \aggre\ and  \imredc\ and proposition \imreds), with the Kac-Moody algebra $\hat\gothg$ (see corollary \rdinsq).

\n Proposition \imreds\ and the structure of the root system of $\hat\gothg$ (see remark \rsy) imply that in order to prove that $\tilde\psi_1\big|_{L_+}$ is injective it is enough to show that for all $\alpha\in Q_{0,+}\setminus\{0\}
$, $r\in\N$ $$dim_{\C}(L_+)_{\alpha+r\delta}\leq dim_{\C}\hat\gothg_{\alpha+r\delta}=\begin{cases}1&{\rm{if}}\ \alpha+r\delta\in\Phi_+^{{\rm{re}}}\cr 0&{\rm{otherwise}}.\end{cases}$$

\n Notice that the results of section \spquno\ imply the following:
\vskip .3 truecm
\prop{\ddd}
\n If $i\in I_0$, $r\in\N$ then: 
$$dim(L_+)_{\alpha_i+r\delta}\leq\begin{cases}1&{\rm{if}}\ \tilde d_i|r\cr 0&{\rm{otherwise}},\end{cases}\leqno{(D1)}$$
$$dim(L_+)_{2\alpha_i+r\delta}\leq\begin{cases}1&{\rm{if}}\ (X_{\tilde n}^{(k)},d_i)=(A_{2n}^{(2)},1)\ {\rm{and}}\ 2\not|r\cr 0&{\rm{otherwise}},\end{cases}\leqno{(D2)}$$
$$dim(L_+)_{h\alpha_i+r\delta}=0\ \ {\rm{if}}\ h>2.\leqno{(D3)}$$
\dim
\n $(L_+)_{\alpha_i+r\delta}=<x_{i,r}^+>$ for all $r\in\N$ and,
for all $h>1$,  $r\in\N$, 
$$(L_+)_{h\alpha_i+r\delta}=<[(L_+)_{(h-1)\alpha_i+r_1\delta},x_{i,r_2}^+]|r_1+r_2=r>;$$ in particular:

\n i) $(L_+)_{h\alpha_i+r\delta}=\{0\}$ for all $r\in\N$ implies $L_{\tilde h\alpha_i+r\delta}=\{0\}$ for all $r\in\N$ $\forall \tilde h\geq h$;

\n ii) $(L_+)_{h\alpha_i+r\delta}=\{0\}$ for all $h>0$, $i\in I_0$ and $r\in\N$ such that $\tilde d_i\not|r$; in particular $(D1)$ holds;

\n iii) $(D2)$ follows from lemma \imleu;

\n iv) $(D3)$ follows from lemma \imled.

\vskip .3 truecm
\n In order to generalize this result to all the roots we embed $L_+$ into a $\gothg_0$-module $L$: this structure provides the symmetries that allow to determine easily the needed dimensions of the homogeneous components of $L_+$.

\vskip .3 truecm
\ddefi{\loab}
\n i) $L_0$ is the abelian Lie algebra generated by $\{h_{i,r}|i\in I_0, r\in\N\}$ with relations $h_{i,r}=0$ if $\tilde d_i\not|r$ (hence $\{h_{i,r}|i\in I_0, \tilde d_i|r\in\N\}$ is a basis of $L_0$);

\n ii) $L_-=(L_+)^{op}$; 

\n iii) for all $i\in I_0$, $r\in\N$, $x_{i,r}^-$ denotes $-x_{i,r}^+$ as an element of $L_-$;

\n iv) $L=L_-\oplus L_0\oplus L_+$;

\n v) $\sigma:L\to L$ is the linear map defined by $L=L_-\oplus L_0\oplus L_+\ni(y,h,x)\mapsto(x,h,y)\in L_-\oplus L_0\oplus L_+=L$ (in particular $\sigma(h_{i,r})=h_{i,r}$, $\sigma(x_{i,r}^{\pm})=-x_{i,r}^{\mp}$).

\vskip .3 truecm
\rem{\lopm}{}
\n i) $L_0=L_0^{op}$ as Lie algebras
(since $L_0$ is abelian); 

\n ii) $\sigma\big|_{L_+}:L_+\to L_-$, $\sigma\big|_{L_-}:L_-\to L_+$ and $\sigma\big|_{L_0}=id_{L_0}:L_0\to L_0$ are anti-isomorphisms of Lie algebras.

\vskip .3 truecm
\rem{\lgemb}{}

\n i) $\gothh_0\ni h_i\mapsto h_{i,0}\in L_0$ defines a homomorphism of Lie algebras; 

\n ii) $\gothg_{0,+}\ni e_i\mapsto x_{i,0}^+\in L_+$ defines a homomorphism of Lie algebras, hence it induces an action of $\gothg_{0,+}$ on $L_+$ ($e_i\mapsto{\rm{ad}}_{L_+}x_{i,0}^+$), and  ${\rm{ad}}_{L_+}x_{i,0}^+$ is locally nilpotent; 

\n iii) $\gothg_{0,-}\ni f_i\mapsto x_{i,0}^-\in L_-$ defines a homomorphism of Lie algebras, hence it induces an action of $\gothg_{0,-}$ on $L_-$ ($f_i\mapsto{\rm{ad}}_{L_-}x_{i,0}^-$), and  $
{\rm{ad}}_{L_-}x_{i,0}^-$
is locally nilpotent.

\vskip .3 truecm

\prop{\lach}
\n $h_{i,r}.x_{j,s}^+=B_{ijr}x_{j,r+s}^+
$  with $B_{ijr}=\begin{cases}0&{\rm{if}}\ \tilde d_{i,j}\not|r\cr 2(2-(-1)^{r})&{\rm{if}}\ (X_{\tilde n}^{(k)},d_i,d_j)=(A_{2n}^{(2)},1,1)\cr a_{ij}&{\rm{otherwise}}\end{cases}$ defines 
a Lie algebra homomorphism $D_+:L_0\to Der(L_+)$.

\n Then $L_0\oplus L_+=L_0\ltimes_{D_+}L_+$ is endowed with a Lie algebra structure (semidirect product of $L_0$ and $L_+$).

\n Of course $(L_0\ltimes_{D_+}L_+)^{op}=L_-\oplus L_0$ is a Lie algebra.
\dim
It is obvious that for all $i\in I_0$, $r\in\N$ the ideal generated by the relations $({\tilde{dr}}_+)$ is stable under the derivation $x_{j,s}^+\mapsto B_{ijr}x_{j,r+s}^+$ (see also \realdrel), hence $h_{i,r}.$ defines a derivation of $L_+$; it is also immediate to see that  $h_{i,r}.=0$ if $\tilde d_i\not |r$ and that $h_{i,r}.h_{j,s}.=h_{j,s}.h_{i,r}.$, hence the map $h_{i,r}\mapsto h_{i,r}.$ induces a Lie algebra homomorphism $D_+:L_0\to Der(L_+)$.

\vskip .3 truecm
\rem{\pmstr}{}
\n i) $\sigma\big|_{L_0\oplus L_+}:L_0\oplus L_+\to L_-\oplus L_0$ and $\sigma\big|_{L_-\oplus L_0}:L_-\oplus L_0\to L_0\oplus L_+$ are anti-isomorphisms of Lie algebras;

\n ii) $\sigma\compo({\rm{ad}}_{L_0\oplus L_+}a)\compo\sigma\big|_{L_-\oplus L_0}=-({\rm{ad}}_{L_-\oplus L_0}\sigma(a))\big|_{L_-\oplus L_0}$ $\forall a\in L_0\oplus L_+$;

\n iii) $\gothh_0\oplus\gothg_{0,+}\to L_0\oplus L_+$ and $\gothg_{0,-}\oplus\gothh_0\to L_-\oplus L_0$ are homomorphisms of Lie algebras (indeed $B_{ij0}=a_{ij}$); in particular they induce actions of $\gothh_0\oplus\gothg_{0,+}$ on $L_0\oplus L_+$ ($h_i\mapsto{\rm{ad}}_{L_0\oplus L_+}h_{i,0}$, $e_i\mapsto{\rm{ad}}_{L_0\oplus L_+}x_{i,0}^+$) and of $\gothg_{0,-}\oplus\gothh_0$ on $L_-\oplus L_0$ ($f_i\mapsto{\rm{ad}}_{L_-\oplus L_0}x_{i,0}^-$, $h_i\mapsto{\rm{ad}}_{L_-\oplus L_0}h_{i,0}$);

\n iv) for all $i\in I_0$ ${\rm{ad}}_{L_0\oplus L_+}x_{i,0}^+\big|_{L_0}$ maps $L_0$ in $L_+$ ($L_+$ is and ideal of $L_0\oplus L_+$), hence 
${\rm{ad}}_{L_0\oplus L_+}x_{i,0}^+$  
is locally nilpotent, since it is locally nilpotent on $L_+$ (see remark \lgemb,ii)); analogously ${\rm{ad}}_{L_-\oplus L_0}x_{i,0}^-$  
is locally nilpotent;

\n v) for all $h\in L_0$ $${\rm {ad}}_{L_0\oplus L_+}h\big|_{L_0}=0={\rm {ad}}_{L_-\oplus L_0}h\big|_{L_0}\ {\rm{and}}\ \sigma\compo {\rm{ad}}_{L_0\oplus L_+}h\compo\sigma\big|_{L_-}=- {\rm{ad}}_{L_-\oplus L_0}h\big|_{L_-};$$ in particular the adjoint actions of $L_0$ on $L_0\oplus L_-$ and on $L_-\oplus L_0$ coincide on $L_0$ and thus define an $L_0$-module structure on $L$ (denoted by $h\mapsto h
_{L}$) such that $\sigma\compo h_{L}\compo\sigma=-h_{L}$; 

\n vi) $\gothh_0(\subseteq L_0)$ acts diagonally on $L$ and trivially on $L_0$; more precisely $L_{\pm}$ and $L_0$, hence $L$, are $Q$-graded ($x_{i,r}^{\pm}\in (L_{\pm})_{\pm\alpha_i+r\delta}=L_{\pm\alpha_i+r\delta}$ and $h_{i,r}\in(L_0)_{r\delta}=L_{r\delta}$) and $h\in\gothh_0$ acts on $L_{\alpha}$ as $\alpha(h)id_{L_{\alpha}}$;

\n vii) the action of $\gothh_0\oplus\gothg_{0,+}$ on $L_0\oplus L_+$ and that of $\gothg_{0,-}\oplus\gothh_0$ on $L_-\oplus L_0$ are obviously homogeneous.

\vskip .3 truecm
\rem{\lpacg}{}

\n We want to provide $L$ with a $\gothg_0$-module structure extending the $\gothh_0$-module structure (remark \pmstr,vi)), compatible with the $\gothh_0\oplus\gothg_{0,\pm}$-module structure on $L_0\oplus L_{\pm}$ (remark \pmstr,iii)), and homogeneous with respect to the $Q$-grading.
\vskip .3 truecm
\rem{\efhom}{}
\n Let $e_{i,L},f_{i,L}:L\to L$ be homogeneous linear maps (that is $e_{i,L}(L_{\alpha})\subseteq L_{\alpha+\alpha_i}$, $f_{i,L}(L_{\alpha})\subseteq L_{\alpha-\alpha_i}$). Then: 

\n i)   the relations $[(h_i)_L,e_{j,L}]=a_{ij}e_{j,L}$, $[(h_i) _L,f_{j,L}]=-a_{ij}f_{j,L}$ are automatically satisfied (because of the diagonal action of $\gothh_0$ on $L$, see remark \pmstr,vi)); 

\n ii) if moreover $e_{i,L}\big|_{L_+}={\rm{ad}}_{L_+}x_{i,0}^+$ $f_{i,L}\big|_{L_-}={\rm{ad}}_{L_-}x_{i,0}^-$, then 
$e_{i,L}$ and 
$f_{i,L}$ 
are locally nilpotent (see remark \lgemb,ii) and iii) and notice that for all $x\in L$ there exists $m\in\N$ such that $e_{i,L}^m(x)\in L_+$, $f_{i,L}^m(x)\in L_-$);

\n iii) 
for all $r\in\N$  $L^{(r)}=\oplus_{\alpha\in Q_0}L_{\alpha+r\delta}$ is 
$e_{i,L}$ and $f_{i,L}$-stable, 
$L=\oplus_{r\in\N}L^{(r)}$.

\vskip .3 truecm
\ddefi{\dvfi}
\n Given $i\in I_0$ let $f_{i,L_+}:L_+\to L_0\oplus L_+$ 
 be the derivation defined on the generators by $f_{i,L_+}(x_{j,r}^+)=-\delta_{i,j}h_{i,r}$ and $e_{i,L_-}:L_-\to L_-\oplus L_0$  be defined by $e_{i,L_-}=\sigma\compo f_{i,L_+}\compo\sigma\big|_{L_-}$. 
\vskip .3 truecm
\prop{\lacf}
\n $f_{i,L_+}$ and $e_{i,L_-}$ are well defined derivations.
\dim{}
\n Obviously if $\rho$ is a relation involving only indices in $I_0\setminus\{i\}$ then $f_{i,L_+}(\rho)=0$; it is also obvious that 
if $\tilde d_j\not|r$ $f_{i,L_+}(x_{j,r}^+)=-\delta_{ij}h_{jr}=0$ (hence $f_{i,L_+}$ preserves relation $(zx)$). 

\n Moreover:

i) if $i\neq j$ $f_{i,L_+}([x_{i,r}^+,x_{j,s}^+])=-a_{ij}x_{j,r+s}^+$ which depends only on $r+s$, hence relation $(x_d)$ is preserved by $f_{i,L_+}$ and symmetrically by $f_{j,L_+}$;

ii) $f_{i,L_+}((-1)^s[x_{i,r}^+,x_{i,s}^+])=(-1)^s(-B_{iir}+B_{iis})x_{i,r+s}^+$ which is zero if $(X_{\tilde n}^{(k)},d_i)\neq (A_{2n}^{(2)},1)$ or $2|r+s$ and in any case depends only on $r+s$, hence relation $(x_{1,2})$ is preserved by $f_{i,L_+}$; 

iii) if $(X_{\tilde n}^{(k)},d_i)=(A_{2n}^{(2)},1)$ then 
$$f_{i,L_+}([[x_{i,r_1}^+,x_{i,r_2}^+],x_{i,r_3}^+])\!\!=\!\![f_{i,L_+}([x_{i,r_1}^+,x_{i,r_2}^+]),x_{i,r_3}^+]+[[x_{i,r_1}^+,x_{i,r_2}^+],f_{i,L_+}(x_{i,r_3}^+)]$$ which, if $2|r_1+r_2$ or $2|r_1+r_2+r_3$, is of course zero by ii) and relation $(x_{1,2})$, while is 
$(-6+2)[x_{i,r_1+r_2}^+,x_{i,r_3}^+]+2[x_{i,r_1+r_3}^+,x_{i,r_2}^+]+2[x_{i,r_1}^+,x_{i,r_2+r_3}^+]=0$ if $2\not| r_1$, $2|r_2$ and $2|r_3$, by relation $(x_{1,2})$.

\n It follows that $f_{i,L_+}$ preserves also relations $(x_{1,2},x_3)$.

\n Furthermore if $k>1$, $a_{ij}=-k$:

iv) $f_{j,L_+}((-1)^{r_2}[[x_{j,s}^+,x_{i,r_1}^+],x_{i,r_2}^+])=-(-1)^{r_2}a_{ji}[x_{i,s+r_1}^+,x_{i,r_2}^+]$ which is zero if $(X_{\tilde n}^{(k)},d_i)\neq(A_{2n}^{(2)},1)$ or $2|s+r_1+r_2$ and depends only on $s+r_1+r_2$ otherwise;

v) $f_{i,L_+}(
[[x_{j,s}^+,x_{i,r_1}^+],x_{i,r_2}^+])=$
$$=
a_{ij}[x_{j,s+r_1}^+,x_{i,r_2}^+]+a_{ij}[x_{j,s+r_2}^+,x_{i,r_1}^+]+B_{iir_2}[x_{j,s}^+,x_{i,r_1+r_2}^+]
$$
(we can suppose, and we are supposing, $\tilde d_j|s$ and $\tilde d_i|r_1,r_2$).

\n Let us distinguish three cases:

\n $k=2$, $X_{\tilde n}^{(k)}\neq A_{2n}^{(2)}$; then $x_{j,s+1}^+=0$, so that
$$f_{i,L_+}(
[[x_{j,s}^+,x_{i,1}^+],x_{i,0}^+])=
-2[x_{j,s}^+,x_{i,1}^+]+2[x_{j,s}^+,x_{i,1}^+]=0
$$
and $f_{i,L_+}(
[[x_{j,s}^+,x_{i,1}^+],x_{i,1}^+]+[[x_{j,s+2}^+,x_{i,0}^+],x_{i,0}^+])=$
$$=2[x_{j,s}^+,x_{i,2}^+]-2[x_{j,s+2}^+,x_{i,0}^+]-2[x_{j,s+2}^+,x_{i,0}^+]+2[x_{j,s+2}^+,x_{i,0}^+]=0;$$
together with i), ii) and iv) this implies the stability of $(x_d,x_{1,2},t_2^{\prime},t_2^{\prime\prime})$ by the action of the $f_l$'s ($l\in I_0$);

\n $X_{\tilde n}^{(k)}=A_{2n}^{(2)}$; then
$f_{i,L_+}(
[[x_{j,s}^+,x_{i,1}^+],x_{i,0}^+]+[[x_{j,s}^+,x_{i,0}^+],x_{i,1}^+])=$
$$=-2[x_{j,s+1}^+,x_{i,0}^+]-2[x_{j,s}^+,x_{i,1}^+]+2[x_{j,s}^+,x_{i,1}^+]+$$
$$-2[x_{j,s}^+,x_{i,1}^+]-2[x_{j,s+1}^+,x_{i,0}^+]+6[x_{j,s}^+,x_{i,1}^+]=0;$$
together with i), ii) and iv) this implies the stability of $(x_d,x_{1,2},s_2)$ by the action of the $f_l$'s;

\n $k=3$; then $f_{i,L_+}(
[[x_{j,s}^+,x_{i,r_1}^+],x_{i,r_2}^+])=$
$$\begin{cases}(-6+2)[x_{j,s}^+,x_{i,r_1+r_2}^+]
&{\rm{if}}\ 3|r_1,\ 3|r_2\cr
(-3+2)[x_{j,s}^+,x_{i,r_1+r_2}^+]&{\rm{if}}\ 3\not|r_1,\ 3|r_2\ {\rm{or}}\ 3|r_1,\ 3\not|r_2\cr
2[x_{j,s}^+,x_{i,r_1+r_2}^+]&{\rm{if}}\ 3\not|r_1,\ 3\not|r_2;\cr
\end{cases}$$
in particular
$$f_{i,L_+}(
[[x_{j,s}^+,x_{i,1}^+],x_{i,1}^+]+2[[x_{j,s}^+,x_{i,2}^+],x_{i,0}^+])=$$
$$=(2-2)[x_{j,s}^+,x_{i,2}^+]=0,$$
$$f_{i,L_+}(
2[[x_{j,s}^+,x_{i,2}^+],x_{i,1}^+]+[[x_{j,s+3}^+,x_{i,0}^+],x_{i,0}^+])=$$
$$=4[x_{j,s}^+,x_{i,3}^+]-4[x_{j,s+3}^+,x_{i,0}^+]=0,$$
$$f_{i,L_+}(
[[x_{j,s}^+,x_{i,2}^+],x_{i,2}^+]+2[[x_{j,s+3}^+,x_{i,1}^+],x_{i,0}^+])=$$
$$2[x_{j,s}^+,x_{i,4}^+]-2[x_{j,s+3}^+,x_{i,1}^+]=0,$$
which, together with i), ii) and iv), implies the stability of $(x_d,x_{1,2},t_3^{\prime},t_3^{\prime\prime},t_3^{\prime\prime\prime})$ by the action of the $f_l$'s;

\n In case $k=3$, $a_{ij}=-3$ 
then:

vi) $f_{j,L_+}(
[[[x_{j,s}^+,x_{i,r_1}^+],x_{i,r_2}^+],x_{i,r_3}^+])=-a_{ji}[[x_{i,s+r_1}^+,x_{i,r_2}^+],x_{i,r_3}^+]=0$;

\n if moreover $3\not|r$ then:

vii)  $f_{i,L_+}(
[[[x_{j,s}^+,x_{i,r}^+],x_{i,0}^+],x_{i,0}^+])=(2a_{ij}+3a_{ii})[[x_{j,s}^+,x_{i,r}^+],x_{i,0}^+]=0$;

\n vi) and vii), together with i) and ii), imply the stability of $(x_d,x_{1,2},x_3,u_3^{\prime},u_3^{\prime\prime})$ by the action of the $f_l$'s.

\n Finally if $i\neq j$ it is well known that 

viii) $f_{l,L_+}(({\rm{ad}}x_{i,0}^+)^{1-a_{ij}}x_{j,s}^+)=0$, which implies the stability of $(serre)$
by the action of the $f_l$'s.

\vskip.3 truecm
\ddefi{\lpacef}
\n Let $e_{i,L},f_{i,L},h_{i,L}:L\to L$ be defined by: 

\n $e_{i,L}\big|_{L_0\oplus L_+}={\rm{ad}}_{L_0\oplus L_+}x_{i,0}^+$ and $f_{i,L}\big|_{L_-\oplus L_0}={\rm{ad}}_{L_-\oplus L_0}x_{i,0}^-$ (see remark \pmstr,iii)); 

\n $f_{i,L}\big|_{L_+}=f_{i,L_+}$ and $e_{i,L}\big|_{L_-}=e_{i,L_-}$ (see definition \dvfi);

\n $h_{i,L}=(h_{i,0})_L$ (see remark \pmstr, v)).
\vskip .3 truecm
\ddefi{\gtg}
\n Define $\tilde\gothg_0$ to be the Lie-algebra generated by $\{e_i,f_i,h_i|i\in I_0\}$ with relations 
$$[h_i,h_j]=0,\ [h_i,e_j]=a_{ij}e_j,\ [h_i,f_j]=-a_{ij}f_j,\ [e_i,f_j]=\delta_{i,j}h_i
\ \ \ \forall i,j\in I_0.$$
\vskip .3 truecm
\lem{\eifi}
\n i) $e_{i,L},f_{i,L}:L\to L$ are homogeneous linear maps, hence  locally nilpotent; 

\n ii) $\sigma\compo f_{i,L}\compo\sigma=e_{i,L}$;

\n iii) $[e_{i,L},f_{j,L}]=\delta_{i,j}h_{i,L}$;

\n iv) $L$ is a $\tilde\gothg_0$-module.
\dim
\n i) follows from remark \efhom,ii), and ii) from remark \pmstr, ii) and from definitions \dvfi\ and \lpacef).

\n iii) By ii) and remark \pmstr, v) it is enough to prove the identity on $L_0\oplus L_+$: by homogeneity $e_{i,L}\compo f_{j,L}$ and $f_{j,L}\compo e_{i,L}$ map $L_0$ in $L_{\alpha_i-\alpha_j}\subseteq L_0$, and in particular in $\{0\}$ if $i\neq j$ while $e_{i,L}\compo f_{i,L}\big|_{L_0}=\sigma\compo e_{i,L}\compo f_{i,L}\compo\sigma\big|_{L_0}=f_{i,L}\compo e_{i,L}\big|_{L_0}$ because $\sigma\big|_{L_0}=id_{L_0}$; hence $[e_{i,L},f_{j,L}]\big|_{L_0}=0=\delta_{i,j}(h_i)_L\big|_{L_0}$; on the other hand, since $f_{j,L}\big|_{L_+}=f_{j,L_+}:L_+\to L_0\oplus L_+$ is a derivation,
$$f_{j,L}\compo e_{i,L}\big|_{L_+}=f_{j,L}\compo {\rm{ad}}_{L_+}x_{i,0}^+=$$
$$=({\rm{ad}}_{L_0\oplus L_+}(f_{j,L}(x_{i,0}^+))+( {\rm{ad}}_{L_0\oplus L_+}x_{i,0}^+)\compo f_{j,L})\big|_{L_+}=$$
$$=(-\delta_{ij}{\rm{ad}}_{L_0\oplus L_+}h_{i,0}+e_{i,L}\compo f_{j,L})\big|_{L_+}=(-\delta_{ij}(h_{i,0})_L+e_{i,L}\compo f_{j,L})\big|_{L_+}=$$
$$=(-\delta_{ij}h_{i,L}+e_{i,L}\compo f_{j,L})\big|_{L_+},$$
which is the claim.

\n iv) is a consequence of iii) together with proposition \pmstr, v) and remark \efhom,i).

\vskip .3 truecm
\lem{\impser}
\n Let $\rho:\tilde\gothg_0\to\gothg\gothl(M
)$ be a $\tilde\gothg_0$-module structure on $M$ 
with weight space decomposition $M=\oplus_{\alpha\in 
\gothh_0^*}M_{\alpha}$ ($\rho(h)\big|_{M_{\alpha}}=\alpha(h)id_{M_{\alpha}}$ for all $h\in\gothh_0$, remarking that $\gothh_0\hookrightarrow\tilde\gothg_0$) and suppose that $\rho(e_i)$, $\rho(f_i)$ are locally nilpotent. 

\n Then $M$ is a $\gothg_0$-module.
\dim
\n Let $i\neq j\in I_0$: we want to prove that $\rho({\rm{ad}}(e_i)^{1-a_{ij}}(e_j))=0$ and $\rho({\rm{ad}}(f_i)^{1-a_{ij}}(f_j))=0$.

\n a) Given $x\in M$ homogeneous, the subspace $M_{x}=<\rho(e_i)^r(x),\rho(f_i)^r(x)|r\in\N>$ is finite dimensional and $e_i,f_i,h_i$-stable; 

\n b) for $\tilde M\subseteq M$ finite dimensional there exists $r\in\N$ such that 
$\rho(e_i)^{r}\big|_{\tilde M}=0$; in particular $\exists r_x\in\N$ such that $\rho(e_i)^{r_x}\big|_{M_x}=0$,
$\rho(e_i)^{r_x}\big|_{\rho(e_j)(M_x)}=0$;

\n c) for $r\in\N$ $\rho({\rm{ad}}(e_i)^r(e_j))=\sum_{u=0}^r{r\choose u}\rho(e_i)^{r-u}\rho(e_j)\rho(e_i)^u$; in particular if $r\geq2r_x-1$ $\rho({\rm{ad}}(e_i)^r(e_j))\big|_{M_x}=0$; 

\n d) for $r\in\N$ $[e_i,{\rm{ad}}(e_i)^r(e_j)]={\rm{ad}}(e_i)^{r+1}(e_j)$ and 
$[f_i,{\rm{ad}}(e_i)^r(e_j)]=-r(a_{ij}+r-1){\rm{ad}}(e_i)^{r-1}(e_j)$;

\n e) let $Y=\{r\in\N|\rho({\rm{ad}}(e_i)^r(e_j))\big|_{M_x}=0\}$; then $2r_x-1\in Y\neq\emptyset$, $r\in Y\Rightarrow r+1\in Y$ and $r\in Y\setminus\{0,1-a_{ij}\}\Rightarrow r-1\in Y$; in particular $1-a_{ij}\in Y$ and $\rho({\rm{ad}}(e_i)^{1-a_{ij}}(e_j))(x)=0$. 

\n Then $\rho({\rm{ad}}(e_i)^{1-a_{ij}}(e_j))=0$.

\n Composing $\rho$ with the Lie-automorphism of $\gothg_0$ defined by $e_i\mapsto -f_i$, $f_i\mapsto -e_i$, $h_i\mapsto -h_i$, we get that also $\rho({\rm{ad}}(f_i)^{1-a_{ij}}(f_j))=0$.
\vskip .3 truecm
\cor{\lrgm}
\n $L$ is a $\gothg_0$-module; $L^{(r)}$ is a $\gothg_0$-module for all $r\in\N$.
\dim
The claim is a straightforward consequence of remark \pmstr, vi), of lemma \eifi, i) and iv) and of lemma \impser.
\vskip .3 truecm
\lem{\genw}
\n Let $\gothg$ be a Lie algebra, $\gothh\subseteq\gothg$ a subalgebra, $M$ a $\gothg
$-module with $M=\oplus_{\alpha\in \gothh
^*}M_{\alpha}$, $M_{\alpha}=\{m\in M|h.m=\alpha(h)m\ \forall h\in\gothh\}$. Let  $\tau\in Aut_{Lie}(\gothg
)$, $\varphi\in GL(M)$ be such that: 

\n i) $\tau(\gothh
)=\gothh
$.

\n ii) $\varphi(y.m)=\tau(y).\varphi(m)$ $\forall y\in\gothg
$, $m\in M$; 

\n Then $\tau.=(\tau\big|_{\gothh
}^{-1})^*\in GL(\gothh
^*)$ and $\varphi(M_{\alpha})=M_{\tau.\alpha}$ for all $\alpha\in\gothh
^*$. 

\n In particular $P_M=\{\alpha\in \gothh
^*|M_{\alpha}\neq\{0\}\}$ is $\tau.$-stable and 
$dim(M_{\alpha})=dim(M_{\tau.\alpha})$ for all $\alpha\in P_M$. 
\vskip .3 truecm
\lem{\gnaut}
\n Let $\gothg$ be a Lie-algebra, $M$ be a $\gothg$-module and 
$x\in\gothg$ be such that ${\rm{ad}}x$ and $x_M$ are nilpotent ($x_M$ denotes the map $m\mapsto x.m$), $\tau^{(x)}=exp({\rm{ad}}x)$, $\varphi^{(x)}=exp(x_M)$. Then $\tau^{(x)}\in Aut_{Lie}(\gothg)$,  $\varphi^{(x)}\in GL(M)$ and 
$\varphi^{(x)}(y.m)=\tau^{(x)}(y).\varphi^{(x)}(m)$ $\forall y\in\gothg
$, $m\in M$.
Moreover if $x_1$,..., $x_r\in\gothg$ are such that ${\rm{ad}}x_i$ and $(x_i)_M$ are nilpotent and we set 
$\tau=\tau^{(x_1)}\compo...\compo\tau^{(x_r)}$, 
$\varphi=\varphi^{(x_1)}\compo...\compo\varphi^{(x_r)}$, we still have 
$\tau\in Aut_{Lie}(\gothg)$, $\varphi\in GL(M)$, $\varphi(y.m)=\tau(y).\varphi(m)$ $\forall y\in\gothg$, $m\in M$.
\dim
It is a straightforward consequence of the well known identity
$$x_M^n(y.m)=\sum_{r=0}^n{n\choose r}({\rm{ad}}x)^r(y).x_M^{n-r}(m).$$
\vskip .3 truecm
\rem{\efexp}{}
\n For all $r\in\N$ let us consider the $\gothg_0$-module $L^{(r)}$ and the elements $e_i$, $f_i\in\gothg_0$. Let 
$$\tau_i=exp({\rm{ad}}e_i)exp(-{\rm{ad}}f_i)exp({\rm{ad}}e_i),$$ $$\varphi_i=exp(e_{i,L}\big|_{L_r})exp(-f_{i,L}\big|_{L_r})exp(e_{i,L}\big|_{L_r}).$$
Then it is well known and obvious (from lemmas \eifi,i), \lrgm\ and \gnaut) that $\tau_i\in Aut_{Lie}(\gothg)$ and $\varphi_i\in GL(M)$ are well defined and $\varphi_i(y.m)=\tau_i(y).\varphi_i(m)\ \forall y\in\gothg_0,\ m\in M$.

\n It is also well known (see \kac) that $\tau_i(\gothh_0)=\gothh_0$ and in fact $\tau_i\big|_{\gothh_0}=s_i\in W_0$,
hence by lemma \genw\ $\{\alpha\in \gothh_0^*|(L^{(r)})_{\alpha}(=L_{\alpha+r\delta})\neq \{0\}\}$ is $W_0$-stable and $dimL_{w(\alpha)+r\delta}=dimL_{\alpha+r\delta}$ for all $\alpha\in \gothh_0^*$, $r\in\N$, $w\in W_0$.

\n Recall that $\{\alpha\in \gothh_0^*|L_{\alpha+r\delta}\neq \{0\}\}\subseteq Q_{0,+}\cup(-Q_{0,+})$.

\vskip .3 truecm
\lem{\porb}
\n Let $P\subseteq Q_{0,+}\cup(-Q_{0,+})$ be $W_0$-stable. Then any $\alpha\in P$ is $W_0$-conjugate to an integer multiple of a simple root.

\dim
Let $\alpha\in P\setminus\{0\}$ and take $\beta\in W_0.\alpha\cap Q_{0,+}(\neq\emptyset$ because there exists $\tilde w\in W_0$ such that $\tilde w(Q_{0,+})=-Q_{0,+})$ of minimal height. Since $(\beta|\beta)>0$ there exists $i\in I_0$ such that $(\beta|\alpha_i)>0$, so that, by the choice of $\beta$, $s_i(\beta)\in-Q_{0,+}$. This implies $\beta$ to be a multiple of $\alpha_i$. 
\vskip .3 truecm
\n Let us now come to our point.
\vskip .3 truecm
\prop{\ult}
\n Given $\alpha\in Q_0$ $dimL_{\alpha+r\delta}\leq\begin{cases}1&{\rm{if}}\ \alpha+r\delta\in\Phi\cr 0&{\rm{otherwise}}.\end{cases}$
\dim
We have already proved (see remark \efexp\ and lemma \porb) that $dimL_{\alpha+r\delta}=0$ if $\alpha\not\in\cup_{h>0}h\Phi_{0
}$. 

\n By 
remark \efexp\ it is then enough to prove the claim when $\alpha$ is an integer multiple of a simple root: but this is nothing but $(D1)$,  $(D2)$,  $(D3)$, see proposition \ddd.

\vskip .3 truecm
\cor{\rdinsq}
\n $\tilde\psi_1\big|_{L_+}:L_+\to(\gothg_+\otimes\C[t])^{\chi}$ is an isomorphism of Lie algebras.
\dim
The claim is a consequence of propositions \imreds\ and \ult.
\vskip .5 truecm
\noindent {\bf{\cncl.\ {\bf CONCLUSIONS.}}}
\vskip .5truecm

\n In this section we point out and underline the several consequences of corollary \rdinsq. They include the main result ($\psi$ is an isomorphism) together with other results which is worth evidentiating, both about the Drinfeld realization of affine quantum algebras and the affine Kac-Moody algebras.

\vskip .3 truecm
\teo{\tttu}
\n $\psi:\udr\to\udj$ is an isomorphism. This means that the affine quantum algebras $\udj$ (Drinfeld and Jimbo presentation) and $\udr$ (Drinfeld realization) are different presentations of the same algebra $\uq=\uq(X_{\tilde n}^{(k)})$, the affine quantum algebra of type $X_{\tilde n}^{(k)}$.

\vskip .3 truecm
\teo{\tttd}
\n The product induces an isomorphism
$$\uq^{Dr,-}\otimes\uq^{Dr,0}\otimes\uq^{Dr,+}\cong\udr=\uq$$
(triangular decomposition of the Drinfeld realization, or Drinfeld triangular decomposition of the affine quantum algebra).
\vskip .3 truecm
\n As remarked above (see remark \repsit) the Drinfeld triangular decomposition is essentially different from the Drinfeld and Jimbo triangular decomposition (remark \qremu, ii)). Their precise connection is described in proposition \tttq.
\vskip .3 truecm

\prop{\tttq}
\n i) $\uq^{Dr,+}\cap\uq^{DJ,+}=\uq^{Dr,+,+}$; it is the $\C(q)$-linear span of the ordered monomials in the $E_{r\delta+\alpha}$'s with $\alpha\in Q_{0,+}$, $r\geq 0$;

\n ii) $\uq^{Dr,0}\cap\uq^{DJ,+}=\uq^{Dr,0,+}$;

\n iii) $\tilde\uq^{Dr,-}\cap\uq^{DJ,+}$ is the $\C(q)$-linear span of the ordered monomials in the $E_{r\delta-\alpha}$'s with $\alpha\in Q_{0,+}$, $r> 0$;

\n iv) $\uq^{DJ,+}\cong(\tilde\uq^{Dr,-}\cap\uq^{DJ,+})\otimes_{\C(q)}\uq^{Dr,0,+}\otimes_{\C(q)}\uq^{Dr,+,+}$;

\n v) $\uq^{Dr,+}\cap\tilde\uq^{DJ,-}$ is the $\C(q)$-linear span of the ordered monomials in the $F_{r\delta-\alpha}K_{r\delta-\alpha}$'s with $\alpha\in Q_{0,+}$, $r>0$;

\n vi) $(X_{i,r}^+|i\in I_0, r<0)=\uq^{Dr,+}\cap\tilde\uq^{DJ,-}$ in case $A_1^{(1)}$ and
$(X_{i,r}^+|i\in I_0, r<0)\subsetneq\uq^{Dr,+}\cap\tilde\uq^{DJ,-}$ otherwise;

\n vii) $\uq^{Dr,+}\cap\tilde\uq^{DJ,-}\subsetneq
(X_{i,r}^+|i\in I_0, r\leq 0)$;

\n viii) 
$\uq^{Dr,+}=\uq^{Dr,+,+}\otimes(\uq^{Dr,+}\cap\tilde\uq^{DJ,-})$.
\dim
With the notations of corollary \pip\ we have $\uq^{DJ,+}=U_-\otimes U_0\otimes U_+$ with
$U_-\subseteq\tilde\uq^{Dr,-}$, $U_0\subseteq\uq^{Dr,0}$, $U_+\subseteq\uq^{Dr,+}$, which implies i), ii), iii), iv). 

\n v) is equivalent to iii).

\n In vi) the inclusions are obvious, as well as the claim in case $A_1^{(1)}$; in the other cases there exist indices $i,j\in I_0$ such that $\delta-(\alpha_i+\alpha_j)$ is a root, then 
$F_{\delta-(\alpha_i+\alpha_j)}\in\uq^{Dr,+}\cap\tilde\uq^{DJ,-}$ while $(X_{i,r}^+|i\in I_0, r<0)_{-\delta+(\alpha_i+\alpha_j)}=(0)$.

\n In vii) the inequality is obvious ($X_{i,0}^+\not\in\uq^{Dr,+}\cap\tilde\uq^{DJ,-}$); for the inclusion, consider the subalgebra of $\uq$ generated by $(\uq^{Dr,+}\cap\tilde\uq^{DJ,-})$ and by the $X_{i,0}^+$'s: it is isomorphic to 
$(\uq^{Dr,+}\cap\tilde\uq^{DJ,-})\otimes(X_{i,0}^+|i\in I_0)$, thanks to the triangular decomposition of $\udj$ and to the fact that for all $x\in(\uq^{Dr,+}\cap\tilde\uq^{DJ,-})_{\alpha},\ i\in I_0$ we have $[X_{i,0}^+,x]_{q^{(\alpha|\alpha_i)}}\in\uq^{Dr,+}\cap\tilde\uq^{DJ,-}$, which implies that 
$(\uq^{Dr,+}\cap\tilde\uq^{DJ,-})\otimes(X_{i,0}^+|i\in I_0)$ is not only a $(\uq^{Dr,+}\cap\tilde\uq^{DJ,-})$-module, but also stable by left multiplication by the $X_{i,0}^+$'s, hence a subalgebra of $\uq$; but
of course it contains $(X_{i,r}^+|i\in I_0, r\leq 0)$, so that in order to prove that $(\uq^{Dr,+}\cap\tilde\uq^{DJ,-})\otimes(X_{i,0}^+|i\in I_0)$ and $(X_{i,r}^+|i\in I_0, r\leq 0)$ are equal (which is the claim) it is enough to compare the dimensions of their homogeneous components: for all $m\geq 0$, $\alpha\in Q_{0,+}$
$${\rm{dim}}((\uq^{Dr,+}\cap\tilde\uq^{DJ,-})\otimes(X_{i,0}^+|i\in I_0))_{\alpha-m\delta}=$$
$$=\sum_{\beta\in Q_{0,+}}{\rm{dim}}(\uq^{Dr,+}\cap\tilde\uq^{DJ,-})_{\alpha-m\delta-\beta}{\rm{dim}}(X_{i,0}^+|i\in I_0)_{\beta}=$$
$$=\#\{(m_1\delta-\gamma_1\preceq...\preceq m_s\delta-\gamma_s)|m_u>0,\gamma_u\in Q_{0,+}, \sum_u m_u=m, \sum_u\gamma_u=\alpha-\beta\}\cdot$$
$$\cdot\#\{(\gamma_1^0\preceq...\preceq\gamma_{\tilde s}^0)|\gamma_u^0\in Q_{0,+}, \sum_u\gamma_u^0=\beta\}=$$
$$=\#\{(m_1\delta+\gamma_1\preceq...\preceq m_s\delta+\gamma_s)|m_u>0,\gamma_u\in Q_{0,+}, \sum_u m_u=m, \sum_u\gamma_u=\alpha-\beta\}\cdot$$
$$\cdot\#\{(\gamma_1^0\preceq...\preceq\gamma_{\tilde s}^0)|\gamma_u^0\in Q_{0,+}, \sum_u\gamma_u^0=\beta\}=$$
$$=\#\{(m_1\delta+\gamma_1\preceq...\preceq m_s\delta+\gamma_s)|m_u\geq 0,\gamma_u\in Q_{0,+}, \sum_u m_u=m, \sum_u\gamma_u=\alpha\}=$$
$$={\rm{dim}}\u_{q,m\delta+\alpha}^{Dr,+,+}={\rm{dim}}(X_{i,r}^+|i\in I_0, r\geq 0)_{m\delta+\alpha}={\rm{dim}}(X_{i,r}^+|i\in I_0, r\leq 0)_{-m\delta+\alpha};$$
this chain of equalities follows from the PBW-bases of $\uq^{Dr,+}\cap\tilde\uq^{DJ,-}$ (see v)), of  $(X_{i,0}^+|i\in I_0)$ and of 
$\uq^{Dr,+,+}$ (see i)) and from the isomorphism between $(X_{i,r}^+|i\in I_0, r\leq 0)$ and $(X_{i,r}^+|i\in I_0, r\geq 0)$ (see \realdrel).

\n viii) Thanks to i) and to remarks \grgr, x) and \drromrem, vii)
the claim follows remarking that for $r\leq s\leq 0$ $\lambda ^N(\beta_r)\in-Q_+$ implies $\lambda ^N(\beta_s)\in-Q_+$ (hence
$T_{\lambda}^N(E_{\beta_r})\in\tilde\uq^{DJ,-}\Rightarrow T_{\lambda}^N(E_{\beta_s})\in\tilde\uq^{DJ,-}$, by remark \grgr, vi)).

\vskip .3 truecm
\teo{\tttc}
\n i) $\uq^{Dr,+}$ is the $\C(q)$-algebra generated by $\{X_{i,r}^+|i\in I_0,r\in\Z\}$ with relations $(ZX^+)$ and $(DR)$.

\n ii) $\uq^{Dr,+,+}$ is the $\C(q)$-algebra generated by $\{X_{i,r}^+|i\in I_0,r\in\N\}$ with relations $(ZX_+^+)$, $(DR_+)$, $(S_+)$, $(U3_+)$.

\n iii)  The $\a$-subalgebra ${\cal U}_{\a}^{Dr,+} $ of $\udr$ generated by $\{X_{i,r}^+|i\in I_0,r\in\Z\}$ 
is the $\a$-algebra generated by $\{X_{i,r}^+|i\in I_0,r\in\Z\}$ with relations $(ZX^+)$ and $(DR)$ and is a free $\a$-module: it is an integer form of $\uq^{Dr,+}$.

\n iv) The $\a$-subalgebra  of $\udr$ generated by $\{X_{i,r}^+|i\in I_0,r\in\N\}$ is the $\a$-algebra generated by $\{X_{i,r}^+|i\in I_0,r\in\N\}$ with relations $(ZX_+^+)$, $(DR_+)$, $(S_+)$, $(U3_+)$ and it is a free $\a$-module: it is an integer form of $\uq^{Dr,+,+}$.
\dim
iv) is true by corollary \pfp\ and clearly implies ii). Of course iii) implies i). 

\n iii) follows from iv): let 
$\f$ be the $\a$-algebra freely generated by $\{X_{i,r}^+|i\in I_0,r\in\Z\}$, $\i$ the ideal of $\f$ defined by the relations $(DR)$, 
$t:\f\to\f$ the $\a$-automorphism defined by $X_{i,r}^+\mapsto X_{i,r+\tilde d_i}^+$, $\bar t$ the $\a$-automorphism induced by $t$ on $\f/\i$ and $f:\f/\i\to\udr$ the natural homomorphism. Consider also the natural homomorphism $j:\f_+/\i_+\to\f/\i$ (see notation \spnot\ and remark \dremd, iii)).

\n Since of course $f\compo j=f_+$ 
and $f\compo\bar t=(t_1\compo...\compo t_n)^{-1}\compo f$, $f$ is injective thanks to corollary \pfp, i) and to the fact that $\f/\i=\cup_{N\in\N}\bar t^{-N}(j(\f_+/\i_+))$.

\n In order to prove that $\f/\i$ is free over $\a$ it is enough to remark that the image of the (injective) homomorphism $\psi\compo f$ is contained in $\ua$ (see remarks \qremd, iii) and \imredu), which is well known to be a free $\a$-module (see remark \qremd, ii)).

\vskip .3 truecm
\teo{\ttto}
\n Let ${\cal U}_{\a}^{Dr,+}$ be as in  theorem \tttc,iii) and let $\uad$, ${\cal U}_{\a}^{Dr,-}$ and $\uado$ be the $\a$-subalgebras of $\udr$ generated respectively by $\{X_{i,r}^{\pm},k_i^{\pm 1}, C^{\pm 1}|i\in I_0,r\in\Z\}$, 
$\{X_{i,r}^{-}|i\in I_0,r\in\Z\}$ and 
$\{H_{i,r},k_i, C, \tilde C|i\in I_0,r\in\Z\}$,
where $H_{i,0}={k_i-k_i^{-1}\over q_i-q_i^{-1}}$ and $\tilde C={C-C^{-1}\over q-q^{-1}}$. Then:

\n i) $\uad={\cal U}_{\a}^{Dr,-}\otimes_{\a}\uado\otimes_{\a}{\cal U}_{\a}^{Dr,+}$;

\n ii) $\uado=\uq^{Dr,0}\cap\uad$ and ${\cal U}_{\a}^{Dr,\pm}=\uq^{Dr,\pm}\cap\uad$;

\n iii) $\u_{\a}^{Dr,*}$ is an integer form of $\uq^{Dr,*}$; this means that it is a free $\a$-module and that $\uq^{Dr,*}=\C(q)\otimes_{\a}\u_{\a}^{Dr,*}$;

\n iv) $\uado$ is the $\a$-algebra generated by $\{H_{i,r},k_i, C, \tilde C|i\in I_0,r\in\Z\}$ with relations $(ZH)$, $(CUK)$, 
$$[\tilde C,x]=0,\ \ [k_i, H_{j,0}]=0,\ \ [H_{i,0},H_{j,0}]=0
\ \ \ (i,j\in I_0),\leqno{(CUH)}$$
$$k_i(k_i-(q_i-q_i^{-1})H_{i,0})=1,\ \ C\Big(C-(q-q^{-1})\tilde C\Big)=1,\leqno{(IQ)}$$
$(KH)$, 
$$[H_{i,0},H_{j,r}]=0\ \ \ (i,j\in I_0,\ r\in\Z),\leqno{(KQH)}$$ and $(HH)$;

\n v) $\uad$  is the $\a$-algebra generated by $\{X_{i,r}^{\pm},H_{i,r},k_i, C, \tilde C|i\in I_0,r\in\Z\}$ with relations $(ZX^{\pm})$, $(CUK)$, $(CUH)$, $(IQ)$, $(KQH)$, $(KX^{\pm})$,
$$H_{i,0}X_{j,r}^{\pm}=X_{j,r}^{\pm}(\pm[a_{ij}]_{q_i}k_i+q_i^{\mp a_{ij}}H_{i,0}),
\leqno{(HQX^{\pm})}$$
$(HH)$, $(HX^{\pm})$, $(XX)$, $(X1_{
const}^{\pm})$, $(X3_{
const}^{\pm})$, $(S_{
const}^{\pm})$;

\n vi) $\uad=\ua$.
\dim
First of all remark that $\uado\subseteq\uad$, since $H_{i,0}=[X_{i,0}^+,X_{i,0}^-]$ $\forall i\in I_0$ and $\tilde C=k_i([X_{i,1}^+,X_{i,-1}^-]-CH_{i,0})$ if $i\in I_0$ is such that $\tilde d_i=d_i=1$. Moreover $\u_{\a}^{Dr,*}\subseteq\uq^{Dr,*}=\C(q)\otimes_{\a}\u_{\a}^{Dr,*}$, so that ${\cal U}_{\a}^{Dr,-}\otimes_{\a}\uado\otimes_{\a}{\cal U}_{\a}^{Dr,+}\to\uad$ is injective, thanks to theorem \tttd.

\n Let $\v$ be the $\a$-algebra generated by $\{X_{i,r}^{\pm},H_{i,r},k_i, C, \tilde C|i\in I_0,r\in\Z\}$ with relations 
$(ZX^{\pm})$, $(CUK)$, $(CUH)$, $(IQ)$, $(KQH)$, $(KX^{\pm})$,
$(HQX^{\pm})$,
$(HH)$, $(HX^{\pm})$, $(XX)$, $(X1_{
const}^{\pm})$, $(X3_{
const}^{\pm})$, $(S_{
const}^{\pm})$
and $\v^+$, $\v^-$, $\v^0$ 
be the $\a$-subalgebras of $\v$ generated respectively by 
$\{X_{i,r}^{+}|i\in I_0,r\in\Z\}$, 
$\{X_{i,r}^{-}|i\in I_0,r\in\Z\}$ and 
$\{H_{i,r},k_i, C, \tilde C|i\in I_0,r\in\Z\}$.

\n Let 
$\tilde\v^0$ be the $\a$-algebra generated by 
$\{H_{i,r},k_i, C, \tilde C|i\in I_0,r\in\Z\}$ with relations  
 $(ZH)$, $(CUK)$, 
$(CUH)$,
$(IQ)$,
$(KH)$, 
$(KQH)$ and $(HH)$
and $\tilde\v^{0,0}=\a[k_i, C, H_{i,0},\tilde C|i\in I_0]/J$
where $J$ is the ideal generated by the relations $(IQ)$.

\n Then:

\n a) $\v^*\to\u_{\a}^{Dr,*}$ is well defined and surjective ($*\in\{\emptyset,0,+.-\}$);

\n b) $\v^-\otimes\v^0\otimes\v^+\to\v$ is surjective, thanks to relations $(CUK)$, $(CUH)$, $(KX^{\pm})$, $(HQX^{\pm})$, $(HX^{\pm})$, $(XX)$, so that  the commutativity of the diagram
$$\xymatrix{& &\v^-\otimes\v^0\otimes\v^+ \ar[r]\ar[d] &\v\ar[d]\ar[r]&0 \\&0\ar[r]
&\u_{\a}^{Dr,-}\otimes\uado\otimes\u_{\a}^{Dr,+}\ar[r] &\uad\ar[d]& \\& & &0 &}$$
implies i); ii) follows from i);

\n c) $\tilde\v^{0,0}\cong\ua\cap\uq^{DJ,0}$ is a free $\a$-module (well known, see remark \qremd, ii) and vi)), so that $\tilde\v^{0,0}\to\uad$ is injective; moreover the (well defined) maps
$$\a[H_{i,r}|i\in I_0,r<0]\otimes\tilde\v^{0,0}\otimes\a[H_{i,r}|i\in I_0,r>0]\to\tilde\v^0\ \ {\rm and}\ \ \tilde\v^0\to\v^0$$
are surjective (thanks to relations $(CUK)$, $(CUH)$, $(KH)$, $(KQH)$ and $(HH)$) and the composition
$$\a[H_{i,r}|i\in I_0,r<0]\otimes\tilde\v^{0,0}\otimes\a[H_{i,r}|i\in I_0,r>0]\to\tilde\v^0\to\v^0\to\uado$$
is injective; in particular $\tilde\v^0\cong\v^0\cong\uado$ so that iv) holds, and $\uado$ is a free $\a$-module; together with i) and with theorem \tttc, iii) this implies also iii).

\n d) $\u_{\a}^{Dr,+}\to\v$ is well defined (remark \dremd\ holds also on $\a$, see \realdrel) with image in $\v^+$, and $\u_{\a}^{Dr,+}\to\v^+$ is obviously surjective: since $$\u_{\a}^{Dr,+}\to\v^+\to\u_{\a}^{Dr,+}$$ is the identity we have that $\u_{\a}^{Dr,+}\cong\v^+$; then the commutativity of the diagram
$$\xymatrix{&\v^-\otimes\v^0\otimes\v^+ \ar[r]\ar[d]_{\cong} &\v\ar[d]\ar[r]&0 \\
&\u_{\a}^{Dr,-}\otimes\uado\otimes\u_{\a}^{Dr,+}\ar[r]^{\,\,\,\,\,\,\,\,\,\,\,\,\,\,\,\,\,\,\,\,\,\,\,\,\,\,\,\,\,\,\cong}&\uad& }$$
implies that $\uad\cong\v$, that is v) holds;

\n e) since $\ua$ is $t_i$-stable for all $i\in I_0$, $X_{i,r}^{\pm}\in\ua$ for all $(i,r)\in \iz$; it is also clear that $\uado\subseteq\ua$, hence $\uad\subseteq\ua$; on the other hand clearly $C$, $K_i$, $E_i\in\uad$ for all $i\in I_0$, $E_0\in\uad$ (see \realdrel) and $F_i\in\uad$ for all $i\in I$ since $\uad$ is $\Omega$-stable; then $\ua\subseteq\uad$ and vi) follows.

\vskip .3 truecm

\teo{\tttt}
\n Consider the affine Kac-Moody algebra $$\hat\gothg=\hat\gothg^-\oplus\hat\gothg^0\oplus\hat\gothg^+=(\gothg_-\otimes\C[t^{\pm 1}])^{\chi}\oplus\big((\gothg_0\otimes\C[t^{\pm 1}])^{\chi}\oplus\C c\big)\oplus(\gothg_+\otimes\C[t^{\pm 1}])^{\chi}.$$
Then:

\n i) $\hat\gothg^+$ is the Lie algebra generated by $\{x_{i,r}^+|i\in I_0,\tilde d_i|r\in\Z\}$ with relations:
$$[x_{i,r}^+,x_{j,s}^+]\ {\rm{depends\ only\ on\ }}r+s\ \ (i\neq j\in I_0\ {\rm {fixed}});\leqno{
}$$
$$
[x_{i,r}^{+},x_{i,s}^{+}]\!=\!\begin{cases}\!0&\!\!\!\!{\rm{if}}\ (X_{\tilde n}^{(k)},d_i)\neq(A_{2n}^{(2)},1)\ {\rm{or}}\ 2|r+s\cr \!(-1)^h[x_{i,s+h+1}^{+},x_{i,s+h}^{+}]&\!\!\!\!
{\rm{if
}}\ r=s+2h+1;
\end{cases}\leqno{
}$$
$$[[x_{i,r_1}^{+},x_{i,r_2}^{+}],x_{i,r_3}^{+}]=0\ \ ((r_1,r_2,r_3)\in\N^3
).\leqno{
}$$
$$[[x_{j,s}^+,x_{i,1}^+],x_{i,0}^+]=0
,\leqno{
}$$
$$[[x_{j,s}^+,x_{i,1}^+],x_{i,1}^+]=-[[x_{j,s+2}^+,x_{i,0}^+],x_{i,0}^+]\leqno{
}$$
$(k=2,\ a_{ij}=-2, X_{\tilde n}^{(k)}\neq A_{2n}^{(2)})$.
$$[[x_{j,s}^+,x_{i,0}^+],x_{i,1}^+]+[[x_{j,s+1}^+,x_{i,0}^+],x_{i,0}^+]=0\ \ 
\ \ (a_{ij}=-2, 
X_{\tilde n}^{(k)}= A_{2n}^{(2)})
;\leqno{
}$$
$$[[x_{j,s}^+,x_{i,1}^+],x_{i,1}^+]=-2[[x_{j,s}^+,x_{i,2}^+],x_{i,0}^+]\ \ (k=3,\ a_{ij}=-3)\leqno{
}$$
$$2[[x_{j,s}^+,x_{i,2}^+],x_{i,1}^+]=-[[x_{j,s+3}^+,x_{i,0}^+],x_{i,0}^+]\ \ (k=3,\ a_{ij}=-3)\leqno{
}$$
$$[[x_{j,s}^+,x_{i,2}^+],x_{i,2}^+]=-2[[x_{j,s+3}^+,x_{i,1}^+],x_{i,0}^+]\ \ (k=3,\ a_{ij}=-3)\leqno{
}$$
$$({\rm {ad}}x_{i,0}^+)^{1-a_{ij}}(x_{j,s}^+)=0\ \ (i\neq j\ {\rm{with}}\ k=1\ {\rm{or}}\ a_{ij}\geq-1).\leqno{
}$$

\n ii) $\hat\gothg^-=(\hat\gothg^+)^{op}$.

\n iii) $\hat\gothg^0$ is the Lie-algebra generated by $\{h_{i,r},c|i\in I_0,\tilde d_i|r\in\Z\}$ with relations
$$[c,h_{i,r}]=0,\ \ \ [h_{i,r},h_{j,s}]=\delta_{r+s,0}{rB_{ijr}\over d_j}c.$$

\n iv) $\hat\gothg$ is the Lie-algebra generated by $\{x_{i,r}^+,x_{i,r}^-,h_{i,r},c|i\in I_0,\tilde d_i|r\in\Z\}$ with relations
$$[c,a]=0\ \ \ \forall a,$$
$$[h_{i,r},h_{j,s}]=\delta_{r+s,0}{rB_{ijr}\over d_j}c,$$
$$[h_{i,r},x_{j,s}^{\pm}]=\pm B_{ijr}x_{j,r+s}^{\pm}
,$$
$$[x_{i,r}^+,x_{j,s}^-]=\delta_{i,j}\Big(h_{i,r+s}+\delta_{r+s,0}{r-s\over 2d_i}c\Big),$$
$$[x_{i,r\pm 1}^{\pm},x_{i,r}^{\pm}]=0\ \ (X_{\tilde n}^{(k)}=A_1^{(1)}),$$
$$[[x_{i,r\pm 1}^{\pm},x_{i,r}^{\pm}],x_{i,r}^{\pm}]=0\ \ (X_{\tilde n}^{(k)}=A_2^{(2)}),$$
$$({\rm{ad}}x_{i,r}^{\pm})^{1-a_{ij}}(x_{j,s}^{\pm})=0\ \ \ (n>1,\ i\neq j\in I_0).$$
(Remark that relations $[h_{i,r},h_{j,s}]=\delta_{r+s,0}{rB_{ijr}\over d_j}c$ depend on the others).
\dim
\n Since $$\u(\hat\gothg)=\ua/(q-1,K_i-1|i\in I)=\uad/(q-1,C-1,k_i-1|i\in I_0)$$ (see remark \qremd, vii) and theorem \ttto, vi)) the claims follow from theorem \tttc, iii), remark \drromrem, ii) and theorem \ttto, iv) and v).

\end{document}